\documentstyle[12pt]{article}
 \input amssym.def
\input amssym
\begin{document}
\def \Z{\Bbb Z}
\def \C{\Bbb C}
\def \R{\Bbb R}
\def \Q{\Bbb Q}
\def \N{\Bbb N}
\def \wt{{\rm wt}}
\def \tr{{\rm tr}}
\def \span{{\rm span}}
\def \Res{{\rm Res}}
\def \Res{{\rm QRes}}
\def \End{{\rm End}}
\def \E{{\rm End}}
\def \Ind {{\rm Ind}}
\def \Irr {{\rm Irr}}
\def \Aut{{\rm Aut}}
\def \Hom{{\rm Hom}}
\def \mod{{\rm mod}}
\def \ann{{\rm Ann}}
\def \<{\langle}
\def \>{\rangle}
\def \t{\tau }
\def \a{\alpha }
\def \e{\epsilon }
\def \l{\lambda }
\def \L{\Lambda }
\def \g{\gamma}
\def \b{\beta }
\def \om{\omega }
\def \o{\omega }
\def \c{\chi}
\def \ch{\chi}
\def \cg{\chi_g}
\def \ag{\alpha_g}
\def \ah{\alpha_h}
\def \ph{\psi_h}
\def \be{\begin{equation}\label}
\def \ee{\end{equation}}
\def \bl{\begin{lem}\label}
\def \el{\end{lem}}
\def \bt{\begin{thm}\label}
\def \et{\end{thm}}
\def \bp{\begin{prop}\label}
\def \ep{\end{prop}}
\def \br{\begin{rem}\label}
\def \er{\end{rem}}
\def \bc{\begin{coro}\label}
\def \ec{\end{coro}}
\def \bd{\begin{de}\label}
\def \ed{\end{de}}
\def \pf{{\bf Proof. }}
\def \voa{{vertex operator algebra}}

\newtheorem{thm}{Theorem}[section]
\newtheorem{prop}[thm]{Proposition}
\newtheorem{coro}[thm]{Corollary}
\newtheorem{conj}[thm]{Conjecture}
\newtheorem{lem}[thm]{Lemma}
\newtheorem{rem}[thm]{Remark}
\newtheorem{de}[thm]{Definition}
\newtheorem{hy}[thm]{Hypothesis}
\makeatletter
\@addtoreset{equation}{section}
\def\theequation{\thesection.\arabic{equation}}
\makeatother
\makeatletter

\baselineskip=16pt

\begin{center}{\Large \bf The theory of physical superselection sectors in
terms of vertex operator algebra language}

\vspace{0.5cm}
Haisheng Li\\
Department of Mathematics,
University of California,
Santa Cruz, CA 95064
\end{center}

{\bf Abstract} We formulate an interpretation of the theory of physical
superselection
sectors in terms of vertex operator algebra language. Using this
formulation we give a construction of simple current from a primary semisimple
 element of weight one. We then prove
that if a rational vertex operator algebra $V$ has a simple current $M$
satisfying certain conditions, then $V\oplus M$ has a natural rational vertex
 operator (super)algebra structure.
Applying our results to a vertex operator algebra associated to an affine Lie
algebra, we construct its simple currents and the extension by a simple
current.
We also present two essentially equivalent constructions for twisted modules
for an inner automorphism from the adjoint module or any untwisted module.

\baselineskip=16pt

\section{Introduction}

This paper is motivated by some physical papers ([FRS], [MSc]) in the theory of
superselection sectors ([DHR], [HK]). Let $A$ be the observable algebra for
a quantum
field theory and let $\pi_{0}$ be the vacuum representation on $H_{0}$ of $A$.
In general, $A$ admits infinitely many inequivalent irreducible modules, so a
criterion is
needed to rule out the physical unrelevant modules.
Let $\psi$ be an endomorphism of $A$. Then we obtain a representation
$\pi_{0}\psi$ on $H_{0}$
of $A$. If $U$ is an element of $A$ with a left inverse $U^{*}$, then we have
an
endomorphism $\psi_{U}$ of $A$ defined by $\psi_{U} (a)=UaU^{*}$ for any
$a\in A$. Consequently, we have a representation  $\pi_{0}\psi_{U}$  on $H_{0}$
of $A$. In the algebraic theory of superselection sectors
(see, e.g., [HK]), the physical superselection sectors consist of each
equivalent class of
irreducible
representation $\pi$ which is equivalent to $\pi_{0}\psi$ for some endomorphism
$\psi$ of $A$.
If $W_{i}=(H_{0}, \pi_{0}\psi_{i})$ ($i=1,2$) are superselection
sectors, an intertwinner from $W_{1}$ to $W_{2}$ is defined to be a
homomorphism of $A$-modules, and $W=(H_{0}, \pi_{0}\psi_{2}\psi_{1})$
is defined to be the tensor product module of $W_{1}$ with $W_{2}$. Such a
tensor product module
is in general reducible, but it is assumed to be decomposable into irreducible
ones.
Furthermore, if $(H_{0}, \pi_{0}\psi_{i})$ ($i=1,2,3$) are three superselection
sectors, then an intertwining operator of type
$\left(\begin{array}{c}W_{3}\\W_{1},W_{2}\end{array}\right)$ is defined to be
an intertwinner or a homomorphism $\phi$ of $A$-modules from
$W=(H_{0}, \pi_{0}\psi_{2}\psi_{1})$ to $W_{3}$ and
fusion rules had also been defined accordingly.

On the other hand,
vertex operator algebras were introduced by mathematician ([B], [FLM]), which
are essentially chiral algebras formulated in [BPZ] in two-dimensional
conformal field theory.
Vertex operator algebras
provide a powerful algebraic tool for studying the general structure of
 conformal field theory. For vertex operator algebra theory, the notion of
modules, intertwining operators and fusion rules have been defined in [FLM] and
[FHL]. Furthermore, the notions of tensor product for modules have been also
developed in [HL0-2] and [Li4]. The initial purpose of this paper
is to interpret the physical superselection theory
in terms of vertex operator algebra language ([B], [FLM], [FHL], [HL0-2],
[Li4]).

Note that if $\sigma$ is an endomorphism of a vertex operator algebra $V$, then
by definition
([B],[FLM], [FHL]), $\sigma$ preserves both the vacuum and the Virasoro element
so that $\sigma$
preserves each homogeneous subspace of $V$. If $V$ is simple,
 i.e., $V$ is an irreducible $V$-module, it follows from Schur lemma that any
nonzero endomorphism
is a scalar. Then the twist of $V$ by $\sigma$ is isomorphic to $V$. Therefore,
it is impossible to
obtain
all irreducible modules by twisting $V$ unless $V$ is holomorphic, i.e.,
any irreducible $V$-module is isomorphic to $V$.
Having known the above fact, we turn to a certain associative algebra.

For any vertex operator algebra $V$, Frenkel and Zhu [FZ] constructed a
topological ${\bf Z}$-graded
associative algebra $U(V)$, which was called the universal enveloping algebra
of $V$.
Roughly speaking, $U(V)$ is  the associative algebra
with identity generated by all $a_{n}$ (linear in $a$),
for $a\in V, n\in {\bf Z}$ with certain defining relations coming from
the Jacobi identity and the Virasoro algebra relations.
Then there is a natural 1-1 correspondence
between the set of equivalence classes of lower truncated ${\bf Z}$-graded
weak $V$-modules and the set of equivalence classes of continuous
\footnote{It was pointed out by C. Dong that this condition is necessary}
lower truncated ${\bf Z}$-graded  $U(V)$-modules. As an infinite-dimensional
associative algebra, $U(V)$
has many nontrivial (continuous) endomorphisms. For instance, it was
essentially proved in [Z] that
${\rm Aut}U(V)$ contains $PSL(2,{\bf C})$ as a subgroup. It is reasonable to
believe that $U(V)$ should
play the role of the observable algebra in the algebraic quantum field theory.

Let $\Delta(z)\in U(V)\{z\}$ satisfy the conditions (2.7)-(2.10). Then it is
proved that for any
$V$-module $(M,Y_{M}(\cdot,z))$, $(\tilde{M},Y_{\tilde{M}}(\cdot,z)):=
(M,Y_{M}(\Delta(z)\cdot,z))$ is a $V$-module. Consequently, $\Delta(z)$
induces an endomorphism $\psi$ of $U(V)$ defined by
\begin{eqnarray}
\psi(Y(a,z))=Y(\Delta(z)a,z)\;\;\;\mbox{for any }a\in V.
\end{eqnarray}
Our first theorem claims that if $\Delta(z)$ is invertible, $\tilde{M}$ is
isomorphic
to a tensor product module of $M$ with $\tilde{V}$ in the sense of [HL0-2] and
[Li4].
This implies that
$\tilde{V}$ is a simple current ([SY1-2], [FG]), {\it i.e.}, the tensor functor
associated to
 $\tilde{V}$
gives a permutation on the set of equivalence classes of irreducible
$V$-modules. Next, we give a construction for such a $\Delta(z)$. Let
$h$ be an semisimple element of a vertex operator algebra $V$ satisfying
\begin{eqnarray}
L(n)h=\delta_{n,0}h, h(n)h=\delta_{n,1}\gamma {\bf 1}\;\;\;\mbox{for any }
n\in {\bf Z}_{+},
\end{eqnarray}
where $\gamma$ is a rational number. Then we construct a $\Delta(h,z)$ (in
Section 2)
satisfying (2.7)-(2.10).
Applying our results to a vertex operator algebra $L(\ell,0)$ associated
 to an affine Lie algebra $\tilde{{\bf g}}$, we prove that if the fundamental
(dominant integral) weight
$\lambda_{i}$ is cominimal [FG], then for any complex number
$\ell\ne -\Omega$ (the dual
Coxeter number), $L(\ell, \ell\lambda_{i})$ is a simple current for
$L(\ell,0)$. If $\ell$ is a positive integer, this result has been
proved in [FG] by calculating four-point functions.

Note that for all known rational vertex operator operator algebras (the
definition is given
in Section 2), there are only finitely many inequivalent irreducible modules
and
all irreducible modules are exactly those which are needed in conformal field
theory.
For instance, it was proved ([DL], [FZ], [Li2]) that for any positive
integer $\ell$, the set of equivalence classes of irreducible
$L(\ell,0)$-modules is exactly the
set of equivalence classes of unitary highest weight
$\tilde{{\bf g}}$-modules of level $\ell$. It is also known ([DMZ], [W]) that
if
$c=1-\frac{6(p-q)^{2}}{pq}$, where $p,q\in \{2,3,\cdots\}$ and
$p$ and $q$ are relatively prime, then the set of equivalence classes of
irreducible
 $L(c,0)$-modules is exactly the
set of equivalence classes of lowest weight Virasoro modules in the minimal
series
given in [BPZ]. Therefore, at least for a rational vertex operator algebra $V$,
each irreducible $V$-module is a superselection sector so that
we do not need a superselection rule to rule out unrelevant modules.
Based this interpretation we conjecture that each irreducible $V$-module is
isomorphic
some twist of $V$ by an endomorphism of $U(V)$.

The famous moonshine module vertex operator algebra $V^{\natural}$ [FLM] is the
first
mathematically rigorous construction of ${\bf Z}_{2}$-orbiford theory, which is
constructed by using
a vertex operator algebra together with an irreducible ${\bf Z}_{2}$-twisted
module. From the construction of $V^{\natural}$ it is an extension of a certain
vertex operator algebra by a simple module.
Recently Huang [Hua] has given a conceptual construction.
To generalize the construction of $V^{\natural}$, one is facing
two problems, which are the existence of twisted modules and the extension of a
vertex
operator algebra by a simple module, respectively. In this paper, we consider a
relatively
simple case for the extension problem.
Let $h\in V$ satisfying the condition (1.2). Then we prove that under certain
conditions,
$V\oplus \tilde{V}$ has a natural vertex operator (super)algebra structure.
Furthermore,
assuming that $V$ is rational, we classify all irreducible modules
 for $V\oplus \tilde{V}$ and prove that $V\oplus \tilde{V}$
is rational. By restricting ourselves to a vertex operator algebra $L(\ell,0)$
associated to
the affine Lie algebra $\tilde{s}\ell_{2}$, our results imply  that
 $L(\ell,0)\oplus L(\ell,{\ell\over 2})$ is a rational vertex
operator algebra if $\ell$ is a positive integral multiple of $4$ and that
 $L(\ell,0)\oplus L(\ell,{\ell\over 2})$ is a rational vertex
operator superalgebra if $\ell$ is a positive odd integral multiple of $2$.
In the notation of $L(\ell,{\ell\over 2})$, ${\ell\over 2}$ is the spin number
(half of
the weight) in terms of physical language.
These results are known to physicists (cf. [MSe]). All results have been
greatly extended
in [DLM].

In the last section, we consider the construction of twisted modules from $V$.
We present two essentially equivalent constructions for twisted
modules (or sectors) for an inner automorphism from $V$ or any (untwisted)
module.
This may shed some light
for the construction of twisted modules (or sectors) in general.

{\large \bf Acknowledgment}
We thank Professors Dong, Lepowsky and Mason for many useful discussions and
Professor
J. Fuchs for correcting a mistake about relation between simple currents and
cominimal weights as in [DLM]. We apologize for missing any related references
because of the lack of our knowledge.

\section{An interpretation of superselection sectors in terms of vertex
operator algebra
language}

A {\em vertex operator algebra} ([B], [FLM]) is a quadruple $(V, Y, {\bf
1},\omega)$, where
$V$ is ${\bf Z}$-graded vector space $V=\oplus_{n\in {\bf Z}}V_{n}$ satisfying
the conditions:
for any $n\in {\bf Z}$, $\dim V_{n}<\infty$ and
$V_{n}=0$ for sufficiently small $n$;
where $Y$ is a linear map
$$Y(\cdot,z)\rightarrow ({\rm End}V)[[z,z^{-1}]]; Y(u,z)=\sum_{n\in {\bf
Z}}u_{n}z^{-n-1}$$
such that for any $u,\;v\in V$,
\begin{eqnarray}
u_nv=0\qquad\mbox{for}\;n\;\mbox{sufficiently large};\end{eqnarray}
and that the {\it Jacobi identity} holds:
\begin{eqnarray}
\hspace{1cm}& & z_0^{-1}\delta\left({z_1-z_2\over
z_0}\right)Y(u,z_1)Y(v,z_2)-z_0^{-1}\delta\left({z_2-z_1\over
-z_0}\right)Y(v,z_2)Y(u,z_1)\nonumber \\
&=& z_2^{-1}\delta\left({z_1-z_0\over z_2}\right)Y(Y(u,z_0)v,z_2);
\end{eqnarray}
where ${\bf 1}$ is a vector of $V$, called the {\em vacuum}, satisfying the
conditions:
\begin{eqnarray}
Y({\bf 1},z)={\rm id}_V;\;\;Y(u,z){\bf 1}\in
V[[z]]\;\;\mbox{and}\;\;\lim_{z\rightarrow 0}Y(u,z){\bf 1}=u;
\end{eqnarray}
and where $\omega$ is a vector of $V$, called the {\em Virasoro element}, such
that if we set
\begin{eqnarray}
L(z)=Y(\omega,z)=\sum_{n\in {\bf Z}}L(n)z^{-n-2},\end{eqnarray}
the following Virasoro algebra relations hold:
\begin{eqnarray}
[L(m),L(n)]=(m-n)L(m+n)+{1\over
12}(m^3-m)\delta_{m+n,0}(\mbox{rank}\:V)\end{eqnarray}
with (rank $V$)$\in {\bf C},\;m,\;n\in {\bf Z}$; and the following relations
hold:
\begin{eqnarray}
{d\over dz}Y(v,z)=Y(L(-1)v,z),\;\;\;L(0)u=({\rm wt}u)u=nu\end{eqnarray}
for $u\in V_{n}$. This completes the definition.

A {\em weak} $V$-module is a vector space $M$ together with a linear map
$Y_{M}(\cdot,z)$ from $V$ to ${\rm End}M[[z,z^{-1}]]$ such that
 $Y_{M}({\bf 1},z)=id_{M}$,
$Y_{M}(L(-1)a,z)={d\over dz}Y_{M}(a,z)$ for any $a\in V$ and a suitably
adjusted Jacobi
identity holds.

A {\em $V$-module} is a weak $V$-module $M$ which is a ${\bf C}$-graded vector
space
$M=\oplus_{h\in {\bf C}}M_{h}$, where $M_{h}$ is
the eigenspace of $L(0)$ on $M$ with eigenvalue $h$, such that for any $h\in
{\bf C}$,
$\dim M_{h}<\infty$ and $M_{n+h}=0$ for sufficiently large integer $n$.

A {\em lower truncated ${\rm Z}$-graded weak $V$-module} is a weak $V$-module
$M$ together
with a ${\bf Z}$-grading $M=\oplus_{n\in {\bf Z}}M(n)$ such that
$M(n)=0$ for sufficiently small integer $n$ and that
$$a_{n}M(m)\subseteq M(m+k-n-1)\;\;\;\mbox{ for }a\in V_{k}, m,n, k\in {\bf
Z}.$$
If any lower truncated ${\bf Z}$-graded weak $V$-module is completely
reducible, we call $V$
{\em rational} [Z].

Let $W_{i}$ $(i=1,2,3)$ be $V$-modules. Then an intertwining operator of type
$\left(\begin{array}{c}W_{3}\\W_{1} W_{2}\end{array}\right)$ is defined [FHL]
to be a linear map
$I(\cdot, z)$ from $W_{1}$ to ${\rm Hom}_{{\bf C}}(W_{2}, W_{3})\{z\}$ such
that
$I(L(-1)u,z)={d\over dz}I(u,z)$ for $u\in W_{1}$ and a suitably adjusted Jacobi
identity
holds.
Denote by
$I\left(\begin{array}{c}W_{3}\\W_{1} W_{2}\end{array}\right)$ the space of all
intertwining
operators of the indicated type. The dimension of this vector space is called
the
{\em fusion rule} of this type. Let $Irr(V)$ be the
set of equivalence classes of irreducible $V$-modules and for any $V$-module
$M$, denote
the equivalence class of $M$ by $[M]$. For $[M_{i}], [M_{j}], [M_{k}]\in
Irr(V)$, denote
the fusion rule by $N_{ij}^{k}$.
Suppose that all fusion rules among irreducible modules are
finite. The fusion algebra or Verlinda algebra is defined to be the algebra
linearly spanned by
$Irr(V)$ with the multiplication:
$$ [M_{i}]\cdot [M_{j}]=\sum_{k\in Irr(V)}N_{ij}^{k}[M_{k}].$$

Let $V$ be a vertex operator algebra and let $U(V)$ be
the universal enveloping algebra of $V$ constructed by Frenkel and Zhu [FZ].
 What we will use about $U(V)$ is an abstract fact that $U(V)=\oplus_{n\in {\bf
Z}}U(V)_{n}$
 is a ${\bf Z}$-graded
topological associative algebra such that there is a natural 1-1 correspondence
between the
set of
equivalence classes of lower truncated ${\bf Z}$-graded continuous
$U(V)$-module and the
set of equivalence classes of lower truncated ${\bf Z}$-graded weak $V$-module.
Since we will not use any details of $U(V)$ we will not recall the
construction and we refer the interested reader for $U(V)$ to [FZ].

An endomorphism $\rho$ of $U(V)$ is said to be {\em restricted} if for any
$a\in V$ and any
integer $k$, there is an integer $m$ such that
$\rho (a_{n})\in \oplus_{j\le k}^{\infty}U(V)_{j}$ for $n>m$. Then it is clear
that for any
$V$-module $M$, the $\rho$-twist of $M$ is a (weak) $V$-module. From the
analysis made in
introduction we conjecture that any irreducible
$V$-module is isomorphic to a twist of the adjoint module $V$ by a restricted
endomorphism of $U(V)$.

Let $\Delta(z)\in U(V)\{z\}$ satisfy the following
conditions:
\begin{eqnarray}
& &\Delta(z)a\in V[z,z^{-1}];\\
& &\Delta(z){\bf 1}={\bf 1};\\
& &[L(-1), \Delta(z)]=- {d\over dz}\Delta(z);\\
& &Y(\Delta(z_{2}+z_{0})a,z_{0})\Delta(z_{2})
=\Delta(z_{2})Y(a,z_{0})\;\;\;\mbox{for any }a\in V.
\end{eqnarray}
Let $(M,Y_{M}(\cdot,z))$ be any $V$-module. Then
for any $a,b\in V$ we have:
\begin{eqnarray}
& &z_{0}^{-1}\delta\left(\frac{z_{1}-z_{2}}{z_{0}}\right)
Y_{M}(\Delta(z_{1})a,z_{1})Y_{M}(\Delta(z_{2})b,z_{2})\nonumber\\
& &-z_{0}^{-1}\delta\left(\frac{-z_{2}+z_{1}}{z_{0}}\right)
Y_{M}(\Delta(z_{2})b,z_{2})Y_{M}(\Delta(z_{1})a,z_{1})\nonumber\\
&=&z_{2}^{-1}\delta\left(\frac{z_{1}-z_{0}}{z_{2}}\right)
Y_{M}(Y(\Delta(z_{1})a,z_{0})\Delta(z_{2})b,z_{2})\nonumber\\
&=&z_{2}^{-1}\delta\left(\frac{z_{1}-z_{0}}{z_{2}}\right)
Y_{M}(Y(\Delta(z_{2}+z_{0})a,z_{0})\Delta(z_{2})b,z_{2})\nonumber\\
&=&z_{2}^{-1}\delta\left(\frac{z_{1}-z_{0}}{z_{2}}\right)
Y_{M}(\Delta(z_{2})Y(a,z_{0})b,z_{2}).
\end{eqnarray}
The conditions (2.8) and (2.9) imply
\begin{eqnarray}
Y_{M}(\Delta(z){\bf 1},z)={\rm id}_{M},\;\;[L(-1), Y_{M}(\Delta(z)a,z)]
={d\over dz}Y_{M}(\Delta(z)a,z)\;\;\;\mbox{for }a\in V.
\end{eqnarray}
Therefore $(\tilde{M},Y_{\tilde{M}}(\cdot,z))=: (M, Y_{M}(\Delta(z)\cdot,z))$
 is a weak $V$-module.
Then each $\Delta(z)$ gives rise to a restricted endomorphism $\psi$ of $U(V)$
such that $\psi (Y(a,z))=Y(\Delta(z)a,z)$ for any $a\in V$.
 Let $G(V)$ be the
set of all $\Delta(z)$ satisfying the conditions (2.7)-(2.10).

{\bf Lemma 2.1.} {\it Let $\Delta(z)\in G(V)$, let $M^{i}$
$(i=1,2,3)$ be three $V$-modules and let
$I(\cdot,z)$ be an intertwining operator of type
$\left(\begin{array}{c}M^{3}\\M^{1} M^{2}\end{array}\right)$. Then
$\tilde{I}(\cdot,z)=I(\Delta(z)\cdot,z)$ is an intertwining operator of type
$\left(\begin{array}{c}\tilde{M}^{3}\\M^{1} \tilde{M}^{2}
\end{array}\right)$.}

{\bf Proof.} The $L(-1)$-derivative property for $\tilde{I}(\cdot,z)$ follows
from the condition (2.9) immediately. For any
$a\in V, u\in M^{1}$ we have:
\begin{eqnarray}
& &z_{0}^{-1}\delta\left(\frac{z_{1}-z_{2}}{z_{0}}\right)
Y_{\tilde{M}^{3}}(a,z_{1})\tilde{I}(u,z_{2})
-z_{0}^{-1}\delta\left(\frac{z_{2}-z_{1}}{-z_{0}}\right)
\tilde{I}(u,z_{2})Y_{\tilde{M}^{2}}(a,z_{1})
\nonumber\\
&=&z_{0}^{-1}\delta\left(\frac{z_{1}-z_{2}}{z_{0}}\right)
Y_{M^{3}}(\Delta(z_{1})a,z_{1})I(\Delta(z_{2})u,z_{2})\nonumber\\
& &-z_{0}^{-1}\delta\left(\frac{z_{2}-z_{1}}{-z_{0}}\right)
I(\Delta(z_{2})u,z_{2})Y_{M^{2}}(\Delta(z_{1})a,z_{1})\nonumber\\
&=&z_{2}^{-1}\delta\left(\frac{z_{1}-z_{0}}{z_{2}}\right)
I(Y_{M^{1}}(\Delta(z_{1})a,z_{0})\Delta(z_{2})u,z_{2})\nonumber\\
&=&z_{2}^{-1}\delta\left(\frac{z_{1}-z_{0}}{z_{2}}\right)
I(\Delta(z_{2})Y_{M^{1}}(a,z_{0})u,z_{2})\nonumber\\
&=&z_{2}^{-1}\delta\left(\frac{z_{1}-z_{0}}{z_{2}}\right)
\tilde{I}(Y_{M^{1}}(a,z_{0})u,z_{2}).
\end{eqnarray}
Then the proof is complete.$\;\;\;\;\;\Box$

{\bf Lemma 2.2.} {\it Let $\Delta(z)\in G(V)$ and let $\psi$ be a
$V$-homomorphism from a $V$-module $W$ to another $V$-module $M$. Then
$\psi$ is also a $V$-homomorphism from the $V$-module $\tilde{W}$ to
the $V$-module $\tilde{M}$.}

{\bf Proof.} For any $a\in V, u\in W$, we have:
\begin{eqnarray}
\psi(Y_{\tilde{W}}(a,z)u)&=&\psi(Y_{W}(\Delta(z)a,z)u)\nonumber\\
&=&Y_{M}(\Delta(z)a,z)\psi(u)\nonumber\\
&=&Y_{\tilde{M}}(a,z)\psi(u).
\end{eqnarray}
Thus $\psi$ is a $V$-homomorphism from $\tilde{W}$ to
$\tilde{M}$.$\;\;\;\;\Box$

{\bf Lemma 2.3.} {\it Let $\Delta_{1}(z), \Delta_{2}(z)\in G(V)$. Then
$\Delta_{1}(z)\Delta_{2}(z)\in G(V)$.}

{\bf Proof.} By assumption we have:
\begin{eqnarray}
& &\Delta_{1}(z)\Delta_{2}(z){\bf 1}=\Delta_{1}(z){\bf 1}={\bf 1},\\
& &[L(-1), \Delta_{1}(z)\Delta_{2}(z)]\nonumber\\
&=&[L(-1), \Delta_{1}(z)]\Delta_{2}(z)
+\Delta_{1}(z)[L(-1), \Delta_{2}(z)]\nonumber\\
&=&-\left({d\over
dz}\Delta_{1}(z)\right)\Delta_{2}(z)+\Delta_{1}(z)\left({d\over
dz}\Delta_{2}(z)\right)\nonumber\\
&=&-{d\over dz}\left(\Delta_{1}(z)\Delta_{2}(z)\right),
\end{eqnarray}
and
\begin{eqnarray}
& &z_{2}^{-1}\delta\left(\frac{z_{1}-z_{0}}{z_{2}}\right)
Y(\Delta_{1}(z_{1})\Delta_{2}(z_{1})a,z_{0})\Delta_{1}(z_{2})\Delta_{2}(z_{2})
\nonumber\\
&=&z_{2}^{-1}\delta\left(\frac{z_{1}-z_{0}}{z_{2}}\right)
\Delta_{1}(z_{2})Y(\Delta_{2}(z_{1})a,z_{0})\Delta_{2}(z_{2})\nonumber\\
&=&z_{2}^{-1}\delta\left(\frac{z_{1}-z_{0}}{z_{2}}\right)
\Delta_{1}(z_{2})\Delta_{2}(z_{2})Y(a,z_{0})
\end{eqnarray}
for any $a\in V$. Thus $\Delta_{1}(z)\Delta_{2}(z)\in G(V)$.$\;\;\;\;\Box$

It is clear that $id_{V}\in G(V)$, so that $G(V)$ is a semigroup.

{\bf Lemma 2.4.} {\it Let $\Delta(z)\in G(V)$ such that
$\Delta(z)$ has an inverse $\Delta^{-1}(z)\in
U(V)[[z,z^{-1}]]$. Then $\Delta^{-1}(z)\in G(V)$.}

{\bf Proof.} First, we have: $\Delta^{-1}(z){\bf 1}={\bf 1}$. Since
$\Delta(z)\Delta^{-1}(z)=1$, we have:
\begin{eqnarray}
0=\left({d\over dz}\Delta(z)\right)\Delta^{-1}(z)+\Delta(z){d\over
dz}\Delta^{-1}(z).
\end{eqnarray}
Then
\begin{eqnarray}
{d\over dz}\Delta^{-1}(z)
&=&-\Delta^{-1}(z)\left({d\over dz}\Delta(z)\right)\Delta^{-1}(z)\nonumber\\
&=&-L(-1)\Delta^{-1}(z)+\Delta^{-1}(z)L(-1)\nonumber\\
&=&-[L(-1), \Delta^{-1}(z)].
\end{eqnarray}
For any $a\in V$, we have:
\begin{eqnarray}
& &z_{2}^{-1}\delta\left(\frac{z_{1}-z_{0}}{z_{2}}\right)
\Delta^{-1}(z_{2})Y(a,z_{0})\nonumber\\
&=&z_{2}^{-1}\delta\left(\frac{z_{1}-z_{0}}{z_{2}}\right)
\Delta^{-1}(z_{2})Y(\Delta(z_{1})\Delta^{-1}(z_{1})a,z_{0})
\nonumber\\
&=&z_{2}^{-1}\delta\left(\frac{z_{1}-z_{0}}{z_{2}}\right)
Y(\Delta^{-1}(z_{1})a,z_{0})\Delta^{-1}(z_{2}).
\end{eqnarray}
Thus $\Delta^{-1}(z)\in G(V)$. $\;\;\;\;\Box$

Denote by $G^{o}(V)$ the subgroup of all
invertible elements of $G(V)$.
Let $H(V)$ be the subgroup of $G(V)$ consisting of each $\Delta(z)$ such
that $(V,Y(\Delta(z)\cdot,z))$ is isomorphic to $(V, Y(\cdot,z))$.
%
%
%

{\bf Conjecture 2.5.} {\it Let $V$ be a vertex operator algebra. Then for
any irreducible $V$-module $(M,Y_{M}(\cdot,z))$, there is a
$\Delta(z)\in G(V)$ such that $(V, Y(\Delta(z)\cdot,z))$ is isomorphic
to $(M,Y_{M}(\cdot,z))$.}

Recall the definition of tensor product for modules for a vertex
operator algebra $V$ from [Li4] (see [HL0-2] for a different version).

{\bf Definition 2.6 [Li4].} Let $M^{1}$ and $M^{2}$ be two $V$-modules. A
 {\it tensor product} for the ordered pair $(M^{1},M^{2})$ is a pair
$(M,F(\cdot,z))$ consisting of a $V$-module $M$ and an intertwining
operator $F(\cdot,z)$ of type
$\left(\begin{array}{c}M\\M^{1} M^{2}\end{array}\right)$ such that
the following universal
property holds: For any $V$-module $W$ and any intertwining
operator $I(\cdot,z)$ of type
$\left(\begin{array}{c}W\\M^{1} M^{2}\end{array}\right)$, there
exists a unique $V$-homomorphism $\psi$ from $M$ to $W$ such that
$I(\cdot,z)=\psi\circ F(\cdot,z)$. (Here $\psi$ extends canonically to a
linear map from $M\{z\}$ to $W\{z\}$.)

{\bf Proposition 2.7 [Li4].} {\it If $(M,F(\cdot,z))$ is a tensor
product for the
ordered pair $(M^{1},M^{2})$ of $V$-modules, then for any $V$-module
$M^{3}$, ${\rm Hom}_{V}(M,M^{3})$ is linearly isomorphic to the space
of intertwining operators of type
$\left(\begin{array}{c}M^{3}\\M^{1} M^{2}\end{array}\right)$.}

Let $M$ be a module for $V$ and let $\Delta(z)\in G(V)$. Set $W=V$,
$Y_{W}(\cdot,z)=Y_{V}(\Delta(z)\cdot,z)$. In physical references,
essentially $(M, Y_{M}(\Delta(z)\cdot,z))$ is defined to be the tensor product
module of $M$ with $W$.

{\bf Conjecture 2.8.} Let $M$ be a module for $V$ and let
$\Delta(z)\in G(V)$. Then $(\tilde{M},\tilde{I}(\cdot,z)):=
(M,Y_{M}(\Delta(z)\cdot,z))$ is a tensor product of $(M,\tilde{V})$ in the
sense of [HL0-2]
 and [Li4], where
$I(\cdot,z)$ is defined by $I(u,z)a=e^{zL(-1)}Y_{M}(a,-z)u$ for
$a\in V, u\in M$.

{\bf Proposition 2.9.} {\it Let $(W,F(\cdot,z))$ be a tensor product
for a pair $(M^{1},M^{2})$ of $V$-modules and let
$\Delta(z)\,\in\, G^{o}(V)$. Then
$(\tilde{W},\tilde{F}(\cdot,z))$ is a tensor product of the pair
$(M^{1}, \tilde{M}^{2})$.}

{\bf Proof.} From Proposition 2.7
we have an intertwining operator
$\tilde{F}(\cdot,z)=F(\Delta(z)\cdot,z)$ of type
$\left(\,\begin{array}{c}\tilde{W}\\M^{1} \tilde{M}^{2}\end{array}\\,\right)$.
Let $M$ be
any $V$-module and let $I(\cdot,z)$ be any intertwining operator of type
$\left(\,\begin{array}{c}M\\M^{1} \tilde{M}^{2}\end{array}\,\right)$. Then
$I(\Delta(z)^{-1}\cdot,z)$ is an intertwining operator of type
$\left(\begin{array}{c}\hat{M}\\M^{1} M^{2}\end{array}\right)$,
where
$(\hat{M},Y_{\hat{M}}(\cdot,z))=(M,Y_{M}(\Delta(z)^{-1}\cdot,z))$. By the
universal
property of $(W,
F(\cdot,z))$, there is a unique $V$-homomorphism $\psi$ from $W$ to $\hat{M}$
such that $\hat{I}(\cdot,z)=\psi\circ \hat{F}(\cdot,z)$. By Lemma 2.2,
$\psi$ is a $V$-homomorphism from $W$ to $M$. Since $\Delta(z)u$ only
involves finitely many terms, we have: $I(\cdot,z)=\psi\circ
F(\cdot,z)$. It is not difficult to check the uniqueness. Then the proof is
complete.$\;\;\;\;\Box$

{\bf Corollary 2.10.} {\it Let $M$ be a $V$-module and let
$\Delta(z)\in G^{o}(V)$. Then $\tilde{M}$ is isomorphic to the tensor
product module of $M$ with $\tilde{V}$.}

{\bf Proof.} By Proposition 5.1.6 in [Li4], $(M,F(\cdot,z))$ is a tensor
product for the pair $(M,V)$, where $F(\cdot,z)$ is the transpose
intertwining operator of $Y_{M}(\cdot,z)$. By Proposition 2.9,
$(\tilde{M},\tilde{F}(\cdot,z))$ is a tensor product for $(M,\tilde{V})$.
$\;\;\;\;\Box$

The following definition is due to physicists (see for example [SY1-2], [FG]).

{\bf Definition 2.11.} Let $V$ be a vertex operator algebra. An
irreducible $V$-module $M$ is called a {\it simple
current} if the tensor functor ``$M\times $'' is a permutation acting
on the set of equivalence classes of irreducible $V$-modules.

By  Corollary 2.10 we have:

{\bf Theorem 2.12.} {\it For any $\Delta(z)\in G^{o}(V)$, $(V,
Y(\Delta(z)\cdot,z))$ is a simple current $V$-module.}

We will apply this result to vertex operator
algebras associated to affine Lie algebras later. Combining Conjectures 2.5 and
2.8
we formulate the following conjecture.

{\bf Conjecture 2.13.} {\it The fusion algebra or the Verlinda algebra for
vertex
operator algebra $V$ is isomorphic to the group algebra of $G(V)/H(V)$ over the
ground
field ${\bf C}$.}

Next, we give a construction of $\Delta(z)$ from a semisimple
primary element of weight one.
Let $V$ be a vertex operator algebra and let $h\in V$
satisfying the following conditions:
\begin{eqnarray}\label{e2.15}
L(n)h=\delta_{n,0}h,\;
h_{n}h=\delta_{n,1}\gamma {\bf 1}\;\mbox{for any }n\in {\bf Z}_{+},
\end{eqnarray}
where $\gamma$ is a fixed rational number.
Furthermore, we assume that $h_{0}$ semisimply acts on $V$ with
integral eigenvalues. From now on, we also freely use $h(n)$ for $h_{n}$.
For any $\alpha\in {\bf Q}$, set
\begin{eqnarray}
E^{\pm}(\alpha h,z)=\exp\left(\sum_{k=1}^{\infty}\frac{\alpha h(\pm k)}{
k}z^{\mp k}\right).
\end{eqnarray}
Then from [LW] we have:
\begin{eqnarray}
E^{+}(\alpha h,z_{1})E^{-}(\beta h,z_{2})=\left(1-{z_{2}\over
z_{1}}\right)^{-\gamma \alpha \beta} E^{-}(\beta h,z_{2})E^{+}(\alpha h,z_{1}).
\end{eqnarray}

Define
\begin{eqnarray}
\Delta(h,z)=z^{h(0)}\exp\left(\sum_{k=1}^{\infty}\frac{h(k)}{-k}
(-z)^{-k}\right)=z^{h(0)}E^{+}(-h,-z)\in U(V)\{z\}.
\end{eqnarray}

{\bf Proposition 2.14.} {\it Let $V$ be a vertex operator algebra and
let $h$ be an element of $V$ satisfying (\ref{e2.15}). Then $\Delta(h,z)\in
G^{0}(V)$.}

Proposition 2.14 was essentially proved in [Li3], but for completeness,
we present the proof. First we prove

{\bf Lemma 2.15.} {\it Let $h\in V$ satisfying (\ref{e2.15}).
 Then we have:}
\begin{eqnarray}
E^{+}(h,z_{1})Y(a,z_{2})E^{+}(-h,z_{1})
=Y(z_{1}^{h(0)}\Delta(-h,z_{1}-z_{2}),z_{2})
\;\;\;\mbox{{\it for} }a\in V.
\end{eqnarray}

{\bf Proof.} For any $a\in V$, using the formula
$\displaystyle{[h(k),
Y(a,z)]=\sum_{i=0}^{\infty}\,\left(\,\begin{array}{c}k\\i\end{array}\,\right)
\,z^{k-i}Y(h(i)a,z)}$
we obtain
\begin{eqnarray}
& &\left[\sum_{k=1}^{\infty}\frac{h(k)}{k}z_{1}^{-k}, Y(a,z_{2})\right]
\nonumber\\
&=&\sum_{k=1}^{\infty}\sum_{i=0}^{\infty}{1\over k}
\left(\begin{array}{c}k\\i\end{array}\right)
z_{1}^{-k}z_{2}^{k-i}Y(h(i)a,z_{2})\nonumber\\
&=&\sum_{k=1}^{\infty}{1\over k}z_{1}^{-k}z_{2}^{k}Y(h(0)a,z_{2})
+\sum_{k=1}^{\infty}\sum_{i=1}^{\infty}{1\over k}
\left(\begin{array}{c}k\\i\end{array}\right)
z_{1}^{-k}z_{2}^{k-i}Y(h(i)a,z_{2})\nonumber\\
&=&-\log \left(1-{z_{2}\over z_{1}}\right)Y(h(0)a,z_{2})\nonumber\\
& &+\sum_{k=0}^{\infty}\sum_{i=1}^{\infty}{1\over k+i}
\left(\begin{array}{c}k+i\\i\end{array}\right)
z_{1}^{-k-i}z_{2}^{k}Y(h(i)a,z_{2})\nonumber\\
&=&-\log \left(1-{z_{2}\over z_{1}}\right)Y(h(0)a,z_{2})\nonumber\\
& &+\sum_{i=1}^{\infty}\sum_{k=0}^{\infty}{1\over i}(-1)^{k}
\left(\begin{array}{c}-i\\k\end{array}\right)
z_{1}^{-k-i}z_{2}^{k}Y(h(i)a,z_{2})\nonumber\\
&=&-\log \left(1-{z_{2}\over z_{1}}\right)Y(h(0)a,z_{2})
+\sum_{i=1}^{\infty}{1\over i}(z_{1}-z_{2})^{-i}Y(h(i)a,z_{2}).
\end{eqnarray}
Then
\begin{eqnarray}
& &E^{-}(h,z_{1})Y(a,z_{2})E^{-}(-h,z_{1})\nonumber\\
&=&Y\left(\left(1-{z_{2}\over
z_{1}}\right)^{-h(0)}E^{+}(h,z_{1}-z_{2})a,z_{2}\right)
\nonumber\\
&=&Y\left(z_{1}^{h(0)}\Delta(-h,z_{1}-z_{2})a,z_{2}\right).
\;\;\;\;\Box
\end{eqnarray}

{\bf Proof of Proposition 2.14.} Since
$[h(0),Y(u,z)]=Y(h(0)u,z)$ for any $u\in V$,
we have:
\begin{eqnarray}
z^{h(0)}Y(u,z)z^{-h(0)}=Y(z^{h(0)}u,z)\;\;\;\;\mbox{for any }u\in V.
\end{eqnarray}
Then it follows from the construction of $\Delta(h,z)$ and Lemma 2.15 that
$\Delta(h,z)$ satisfies (2.10). Since
$$[L(-1), h(0)]=0,\; [L(-1), h(k)]=-kh(k-1) \;\;\;\mbox{ for }k\in {\bf Z},$$
we obtain
$$[L(-1), \Delta(h,z)]=\sum_{k=1}^{\infty}h(k-1)(-z)^{-k}\Delta(h,z)
={d\over dz}\Delta(h,z).$$
It is clear that $\Delta(h,z)$ satisfies (2.7) and (2.8). Thus
$\Delta(h,z)\in G^{o}(V)$.
$\;\;\;\;\Box$

At the end of this section, we apply our results to some
concrete examples. Let $L$ be a positive-definite even lattice, let $P$ be the
dual
lattice of $L$  and let $V_{L}$ be the vertex operator algebra constructed by
Borcherds [B], Frenkel, Lepowsky and Meurman [FLM]. Then there is a 1-1
correspondence between the set of equivalence classes of irreducible
modules for $V_{L}$ and the set of cosets of $P/L$ ([B], [FLM], [D1]). More
specifically, $V_{P}$ is a $V_{L}$-module with the following
decomposition into irreducible modules:
\begin{eqnarray}
V_{P}=V_{L}\oplus V_{L+\beta_{1}}\oplus \cdots V_{L+\beta_{k-1}}
\end{eqnarray}
where $k=|P/L|$.

{\bf Proposition 2.16.}  {\it Let $\beta\in P$. Then as a $V_{L}$-module,
$(V_{L},Y(\Delta(\beta,z)\cdot,z))$ is isomorphic to the $V_{L}$-module
$V_{L+\beta}$.}

{\bf Proof.} For any $h'\in H={\bf C}\otimes_{{\bf Z}}L$, we have
\begin{eqnarray}
\Delta(\beta,z)h'=\Delta(\beta,z)h'(-1){\bf 1}=h`+z^{-1}\beta(h').
\end{eqnarray}
Then $\tilde{Y}(h',z)=Y(h',z)+z^{-1}\beta(h')$.
Thus $\tilde{V}_{L}$ is a completely reducible module which is isomorphic to
$V_{L}$,
for the Heisenberg algebra $\tilde{H}$ and the set of eigenfunctions of $H(0)$
on
$\tilde{V}_{L}$ is $\beta+L$. Then it follows from [D1] that
$(V_{L},Y(\Delta(\beta,z)\cdot,z))$ is isomorphic to
$V_{L+\beta}$. $\;\;\;\;\Box$

It follows from Proposition 2.16 that all irreducible $V_{L}$-modules
can be obtained by using some $\Delta(\beta,z)$
and that $(V_{L},Y(\Delta(\beta,z)\cdot,z))$ is isomorphic to
$(V_{L},Y(\cdot,z))$ if and only if $\beta\in L$.
It is clear that $\Delta(\beta,z)$ is invertible so that each
irreducible module is a simple current. It is also clear that
$\Delta(\alpha,z)\Delta(\beta,z)=\Delta(\alpha+\beta,z)$ for $\alpha, \beta\in
P$.
Let $\beta_{i}\in P$ $(i=1,2)$. Then $Y(\Delta(\beta_{2},z)\cdot,z)$ is a
nonzero intertwining operator of type
$$\left(\begin{array}{c}(V_{L},Y(\Delta(\beta_{1}+\beta_{2},z)\cdot,z))\\
(V_{L},Y(\Delta(\beta_{1},z)\cdot,z)) (V_{L},
Y(\Delta(\beta_{2},z)\cdot,z))\end{array}\right).$$
Since each irreducible $V_{L}$-module is a simple current, all fusion rules
are either zero or $1$. This result on fusion rules has been obtained in [DL]
using a different method. It is clear that Conjectures 3.5, 3.8 and 3.13 hold
for $V=V_{L}$.


Let ${\bf g}$ be a finite-dimensional simple Lie algebra with a fixed
Cartan subalgebra $H$, let $\{\alpha_{1},\cdots, \alpha_{n}\}$ be a set of
positive roots and
let $\{e_{i}, f_{i}, h_{i}|i=1,\cdots, n\}$
be the Chevalley generators. Let
$\displaystyle{\theta=\sum_{i=1}^{n}a_{i}\alpha_{i}}$ be the highest
positive root and let $\Omega$ be the dual Coxeter number of ${\bf g}$.
Let $\<\cdot,\cdot\>$ be the normalized Killing form on ${\bf g}$
such that $\<\theta, \theta\>=2$. Let $\lambda_{i}$ $(i=1,\cdots,n)$ be the
fundamental
weights of ${\bf g}$ and let $P_{+}$ be the set of dominant integral weights of
${\bf g}$.
Recall from [Hum] that a dominant integral weight $\lambda$ is said to be {\em
minimal} if
it is minimal in $P_{+}$. $\lambda$ is said to be {\em cominimal} [FG] if
$\lambda^{\lor}$ is minimal for the dual Lie algebra.
{}From the table in [K], $\lambda_{i}$ is cominimal if and only if
$a_{i}=1$.
Let $\tilde{{\bf g}}$ be the affine Lie algebra [K]. For any complex number
$\ell$ and
 any weight $\lambda$ of ${\bf g}$,
let $L(\ell, \lambda)$ be the irreducible highest weight $\tilde{{\bf
g}}$-module of
level $\ell$ with lowest weight $\lambda$.
It has been well known (cf. [FZ]) that
$L(\ell,0)$ has a natural vertex operator algebra structure if $\ell\ne
-\Omega$.

{\bf Proposition 2.17.} {\it For any complex number $\ell\ne -\Omega$,
$L(\ell,\ell\lambda_{i})$ is a simple current for $L(\ell,0)$ if $\lambda_{i}$
is
cominimal.}

{\bf Proof.}  Choose $h\in H$ such that
$\alpha_{j}(h)=\langle h,h_{j}\rangle=\delta_{i,j}$ for $1\le j\le n$. Then we
are going
to show that $(V, Y(\Delta(h,z)\cdot,z))$ is isomorphic to
$L(\ell\Lambda_{i})$.
Since $a_{i}=1$, $\theta (h)=a_{i}=1$. By definition we have:
\begin{eqnarray}
& &\Delta(h,z)h_{j}=h_{j}+\ell \delta_{i,j}z^{-1},\; \Delta(h,z)e_{i}=ze_{i},\;
\Delta(h,z)f_{i}=z^{-1}f_{i},\\
& &\Delta(h,z)e_{j}=e_{j},\;
\Delta(h,z)f_{j}=f_{j},\;\Delta(h,z)f_{\theta}=z^{-1}f_{\theta}
\;\;\mbox{ for }j\ne i.
\end{eqnarray}
In other words, the corresponding automorphism $\psi$ of
$U(\tilde{{\bf g}})$ or $U(L(\ell,0))$ satisfying the following conditions:
\begin{eqnarray}
& &\psi(h_{i}(n))=h_{i}(n)+\delta_{n,0}\ell,\;\psi(e_{i}(n))=e_{i}(n+1),\;
\psi(f_{i}(n))=f_{i}(n-1);\\
& &\psi(h_{j}(n))=h_{j}(n),\;\psi(e_{j}(n))=e_{j}(n),\;
\psi(f_{j}(n))=f_{j}(n)
\;\;\;\mbox{for }j\ne i, n\in {\bf Z},\;\;\;\;\;
\end{eqnarray}
and
\begin{eqnarray}
\psi (f_{\theta}(n))=f_{\theta}(n-1)\;\;\;\mbox{for }n\in {\bf Z}.
\end{eqnarray}
Then the vacuum vector ${\bf 1}$ in $(V, Y(\Delta(h,z)\cdot,z))$ is a
highest weight vector of weight $\ell\lambda_{i}$. Thus $(V,
Y(\Delta(h,z)\cdot,z))$ is isomorphic to $L(\ell,\ell\lambda_{i})$ as a
$\tilde{{\bf g}}$-module. By Theorem 2.12, $L(\ell,\ell\lambda_{i})$ is a
simple current. $\;\;\;\;\Box$

{\bf Remark 2.18:} Proposition 2.17 has been proved  in [FG] by calculating
the four point functions and it has also been proved in [F] that those are
all simple currents except for $E_{8}$.

\section{Extension of certain vertex operator algebras}
The famous moonshine moonshine
vertex operator algebra $V^{\natural}$ [FLM] was built as an extension of a
vertex operator
algebra by an irreducible module. (Another conceptual proof can be found in
[Hua].)
One important problem is to determine whether one can have
an extension.
The study of extension of a vertex operator algebra by a self-dual irreducible
module with integral weights was initiated in [FHL] where they proved the
duality or the Jacobi
identity for two module elements and one algebra element.
This section is devoted to the study of extension of
certain vertex operator algebras by a simple current.

Let $W_{i}$ $(i=1,2,3)$ be three $V$-modules and let $I(\cdot,z)$ be
an intertwining operator of type $\left(\begin{array}{c}W_{3}\\W_{1}
W_{2}\end{array}\right)$. The contragredient intertwining operator
$I'(\cdot,z)$ of $I(\cdot,z)$ is an intertwining operator of type type
$\left(\begin{array}{c}W_{2}'\\W_{1}
W_{3}'\end{array}\right)$ defined ([FHL], [HL0-2]) as follows:
\begin{eqnarray}
\langle u_{2}, I'(u_{1},z)u_{3}'\rangle _{2}=
\langle I(e^{zL(1)}e^{\pi iL(0)}z^{-2L(0)}u_{1},z^{-1})u_{2},u_{3}'\rangle_{3}
\end{eqnarray}
for $u_{i}\in W_{1}, u_{i}'\in W_{i}'$. The transpose intertwining
operator $I^{t}(\cdot,z)$ of $I(\cdot,z)$ is an intertwining operator
of type $\left(\begin{array}{c}W_{3}\\W_{2} W_{1}\end{array}\right)$
defined ([FHL], [HL0-2]) as follows:
\begin{eqnarray}
 I^{t}(u_{2},z)u_{1} =e^{zL(-1)}I(u_{1}, e^{\pi i}z)u_{2}\;\;\;\mbox{
for }u_{i}\in W_{i}.
\end{eqnarray}
It was proved ([HL0-2], [Li4]) that
\begin{eqnarray}
I\left(\begin{array}{c}W_{3}\\W_{1} W_{2}\end{array}\right)
\simeq I\left(\begin{array}{c}W_{3}\\W_{2} W_{1}\end{array}\right)
\simeq I\left(\begin{array}{c}W_{2}'\\W_{1} W_{3}'\end{array}\right).
\end{eqnarray}

Let $V$ be a self-dual vertex operator algebra, {\it i.e.}, $V\simeq V'$ as
a $V$-module and let $M$ be an irreducible $V$-module. From [Li1] we have:
\begin{eqnarray}
I\left(\begin{array}{c}M'\\V M\end{array}\right)\simeq {\rm Hom}_{V}(M,M'),\;
\dim {\rm Hom}_{V}(M,M')\le 1.
\end{eqnarray}
Then $\dim I\left(\begin{array}{c}V\\M M\end{array}\right)\le 1$.
Thus $\dim I\left(\begin{array}{c}V\\M M\end{array}\right)=
1$ is equivalent to that $M$ is self-dual. From [Li1], any invariant bilinear
form on $M$ is either symmetric or skew-symmetric. Let $Y(\cdot,z)$ be
a nonzero intertwining operator of type
$\left(\begin{array}{c}V\\M M\end{array}\right)$ and let $(\cdot,\cdot)$ be the
corresponding
bilinear form on $M$. Then we have the
following minor generalization of Proposition 5.6.1 in [FHL]:

{\bf Proposition 3.1.} {\it The bilinear form $(\cdot,\cdot)_{M}$ is
symmetric (resp. skew-symmetric) if and only if}
\begin{eqnarray}
Y(u,z)v=\pm e^{zL(-1)}Y(e^{2\pi i L(0)}v,e^{\pi i}z)u\;\;\;\mbox{ {\it
for }}u,v\in M,
\end{eqnarray}
{\it where ``$+$'' (resp. ``$-$'') corresponds to a symmetric (resp.
skew-symmetric) bilinear form.}

{\bf Proposition 3.2 [FHL].} {\it Let $V$ be a self-dual vertex operator
algebra and let $M$ be a self-dual $V$-module with integral
weights. Then the duality or the Jacobi identity for three elements in
$V\cup M$ consisting of at most two module elements holds.}

Let $V$ be a vertex operator algebra and let
$M$ be a simple current with integral weights. The question is whether
$V\oplus M$ always has a natural vertex operator
algebra structure which extends $V$.
In this section, our mail goal is to prove that if $h\in V$ satisfies
(\ref{e2.15}) such that
$\gamma$ is an integer and $\tilde{\tilde{V}}$ is isomorphic to the
adjoint module $V$, then $V\oplus \tilde{V}$ is a vertex operator
(super)algebra.

{}{\bf From now on, we assume that $h\in V$ satisfies (\ref{e2.15}) such that
$\gamma$ is an integer and $\tilde{\tilde{V}}$ is
 isomorphic to the adjoint module $V$}. From Section 2, for any $V$-module
$(M,Y_{M}(\cdot,z))$, we have a $V$-module
$(\tilde{M}, Y_{\tilde{M}}(\cdot,z))=(M,Y_{M}(\Delta(h,z)\cdot,z))$. This
yields
a linear isomorphism (the identity map) $\psi_{M}$ from $\tilde{M}$
 onto $M$ such that
\begin{eqnarray}\label{e3.6}
\psi_{M}(Y(a,z)u)=Y(\Delta(h,z)a,z)\psi_{M}(u)\;\;\;\mbox{for }a\in V,
u\in \tilde{M}.
\end{eqnarray}
Since $\Delta(h,z)^{-1}=\Delta(-h,z)$, (\ref{e3.6}) is equivalent to
\begin{eqnarray}\label{e3.7}
\psi_{M}\left(Y(\Delta(-h,z)a,z)u\right)=Y(a,z)\psi_{M} (u).
\end{eqnarray}
Let $\pi_{0}$ be an isomorphism from $V$ onto $\tilde{\tilde{V}}$ and set
$\phi_{V}=\psi_{\tilde{V}}\pi_{0}$. Since $\psi_{\tilde{V}}$ is a linear
isomorphism from
$\tilde{\tilde{V}}$ onto $\tilde{V}$, $\phi_{V}$ is a  linear
isomorphism from $V$ onto $\tilde{V}$ such that
\begin{eqnarray}\label{e3.8}
\phi_{V}\left( Y(a,z)b\right)=Y(\Delta(h,z)a,z)\phi_{V}(b),
\end{eqnarray}
or equivalently
\begin{eqnarray}
\phi_{V} (Y(\Delta(-h,z)a,z)b)=Y(a,z)\phi_{V} (b)\;\;\;\mbox{for }a,b\in V.
\end{eqnarray}

Set $\bar{V}=V\oplus \tilde{V}$. For any $u\in \tilde{V}, a\in V$ we define
\begin{eqnarray}
\bar{Y}(u,z)a:=e^{zL(-1)}Y(a,-z)u,
\end{eqnarray}
and for any $u,v\in \tilde{V}$, we define
\begin{eqnarray}
\bar{Y}(u,z)v:=\pi_{0}^{-1}\psi_{\tilde{V}}^{-1}\bar{Y}(\Delta(h,z)u,z)\psi_{V}
(v)=
\phi_{V}^{-1} \bar{Y}(\Delta(h,z)u,z)\psi_{V} (v).
\end{eqnarray}
By Lemma 2.1  $\bar{Y}(\cdot,z)$ are  intertwining operators of types
$\left(\begin{array}{c}\tilde{V}\\ \tilde{V} V\end{array}\right)$ and
$\left(\begin{array}{c}V\\\tilde{V} \tilde{V}\end{array}\right)$, respectively.

{\bf Lemma 3.3.} {\it The following identities hold}
\begin{eqnarray}
& &e^{z(L(1)-h(1))}e^{-zL(1)}=\exp\left(\sum_{k=1}^{\infty}
\frac{h(k)}{k}(-z)^{k}\right)=E^{+}(h,-z),\\
& &e^{z(L(-1)+h(-1))}e^{-zL(-1)}=\exp\left(\sum_{k=1}^{\infty}
\frac{h(-k)}{k}z^{k}\right)=E^{-}(h,z).
\end{eqnarray}

{\bf Proof.} Set
\begin{eqnarray}
X(z)=e^{z(L(1)-h(1))}e^{-zL(1)}\exp\left(\sum_{k=1}^{\infty}
\frac{h(k)}{-k}(-z)^{k}\right).
\end{eqnarray}
Since $X(z)=1$ when $z=0$,
it is sufficient to prove $\displaystyle{{d\over dz}X(z)=0}$.
It follows from the Jacobi identity and (\ref{e2.15}) that $[L(1),
h(k)]=-kh(k+1)$
for $k\in {\bf Z}$. Then
using induction on $n$ we get
\begin{eqnarray}
h(1)L(1)^{n}=\sum_{k=0}^{n}\frac{n!}{(n-k)!}L(1)^{n-k}h(k+1).
\end{eqnarray}
Thus
\begin{eqnarray}\label{e3.12}
h(1)e^{zL(1)}=e^{zL(1)}\sum_{k=1}^{\infty}h(k)z^{k-1}.
\end{eqnarray}
Using (\ref{e3.12}) we obtain
\begin{eqnarray}
& &{d\over dz}X(z)\nonumber\\
&=&e^{z(L(1)-h(1))}(L(1)-h(1))e^{-zL(1)}\exp\left(\sum_{k=1}^{\infty}
\frac{h(k)}{-k}(-z)^{k}\right)\nonumber\\
& &-e^{z(L(1)-h(1))}L(1)e^{-zL(1)}\exp\left(\sum_{k=1}^{\infty}
\frac{h(k)}{-k}(-z)^{k}\right)
+X(z)\sum_{k=1}^{\infty}h(k)(-z)^{k-1}\nonumber\\
&=&-e^{z(L(1)-h(1))}h(1)e^{-zL(1)}\exp\left(\sum_{k=1}^{\infty}
\frac{h(k)}{-k}(-z)^{k}\right)
+X(z)\sum_{k=1}^{\infty}h(k)(-z)^{k-1}\nonumber\\
&=&-e^{z(L(1)-h(1))}e^{-zL(1)}\left(\sum_{k=1}^{\infty}h(k)(-z)^{k-1}\right)
\exp\left(\sum_{k=1}^{\infty}
\frac{h(k)}{-k}(-z)^{k}\right)\nonumber\\
& &+X(z)\sum_{k=1}^{\infty}h(k)(-z)^{k-1}\nonumber\\
&=&0.
\end{eqnarray}
This proves the first identity. Using the symmetry: $L(1)\mapsto
-L(-1), h(k)\mapsto h(-k)$ for $k=1,2,\cdots$, one can obtain the
second identity. $\;\;\;\;\Box$

Since $\Delta(h,z)h=h+z^{-1}\gamma$, we get
\begin{eqnarray}
\psi_{W} h(n)=(h(n)+\delta_{n,0}\gamma)\psi_{W},\;
 \phi_{W}h(n)=(h(n)+\delta_{n,0}\gamma)\phi_{W}
\;\;\;\mbox{ for }n\in {\bf Z}.
\end{eqnarray}
Then from the definition of $E^{\pm}(h,z)$ and $\Delta(h,z)$ we obtain
\begin{eqnarray}
& &\phi_{W} E^{\pm}(h,z)=E^{\pm}(h,z)\phi_{W},\;
\psi_{W} E^{\pm}(h,z)=E^{\pm}(h,z)\psi_{W};\\
& &\psi_{W} \Delta(h,z)=z^{\gamma}\Delta(h,z) \psi_{W},\;
\phi_{W}\Delta(h,z)=z^{\gamma}\Delta(h,z) \phi_{W}.
\end{eqnarray}

{\bf Lemma 3.4.} {\it Let $W$ be any $V$-module. Then the  following identities
hold}
\begin{eqnarray}
& &\psi_{W} e^{zL(-1)}\psi_{W}^{-1}e^{-zL(-1)}=E^{-}(h,z),\;
 e^{zL(-1)}\psi_{W}e^{-zL(-1)}\psi_{W}^{-1}=E^{-}(-h,z),\\
& &\phi_{W} e^{zL(-1)}\phi_{W}^{-1}e^{-zL(-1)}=E^{-}(h,z),\;
 e^{zL(-1)}\phi_{W}e^{-zL(-1)}\phi_{W}^{-1}=E^{-}(-h,z).
\end{eqnarray}

{\bf Proof.} By a simple calculation we get
\begin{eqnarray}
\Delta(h,z)\omega =\omega +z^{-1}h+{1\over 2}\gamma z^{-2},\;
\Delta(-h,z)\omega =\omega -z^{-1}h+{1\over 2}\gamma z^{-2}.
\end{eqnarray}
Then from (\ref{e3.6}) and (\ref{e3.7}) we obtain
\begin{eqnarray}
\psi_{W} L(-1)=(L(-1)+h(-1))\psi_{W},\; L(-1)\psi_{W}=\psi_{W}(L(-1)-h(-1)).
\end{eqnarray}
Thus
\begin{eqnarray}
\psi_{W} e^{zL(-1)}\psi_{W}^{-1}=e^{z(L(-1)+h(-1))},\;
\psi_{W}^{-1} e^{zL(-1)}\psi_{W}=e^{z(L(-1)-h(-1))}.
\end{eqnarray}
By Lemma 3.3 we obtain
$$\psi_{W} e^{zL(-1)}\psi_{W}^{-1}e^{-zL(-1)}=E^{-}(h,z),\;
\phi_{W}^{-1} e^{zL(-1)}\phi_{W}e^{-zL(-1)}=E^{-}(-h,z).$$
Similarly, one can prove the identities for $\phi_{W}$. $\;\;\;\;\Box$

{\bf Lemma 3.5.} {\it (a) For any $a,b\in V$, we have}
\begin{eqnarray}
Y(\psi_{V} \phi_{V}a,z)b
=E^{-}(-h,z)^{2}\psi_{V} \phi_{V}Y(a,z)\Delta(-h,-z)^{2}b.
\end{eqnarray}

{\it (b) For any $u\in \tilde{V}, a\in V$, we have}
\begin{eqnarray}
\bar{Y}(u,z)a=E^{-}(-h,z)Y(\Delta(h,z)\phi_{V}^{-1}
(u),z)\phi_{V}\Delta(-h,-z)a;
\end{eqnarray}

{\it (c) For any $u,v\in \tilde{V}$, we have}
\begin{eqnarray}
\bar{Y}(u,z)v=z^{-\gamma}E^{-}(-h,z)Y(\Delta(h,z)\phi_{V}^{-1}
(u),z)\Delta(-h,-z)\psi_{V}(v).
\end{eqnarray}

{\bf Proof.} Using the skew-symmetry ([B], [FHL]) $Y(u,z)v=e^{zL(-1)}Y(v,-z)u$
and Lemma 3.4
we obtain
\begin{eqnarray}
\bar{Y}(\psi_{V} \phi_{V}a,z)b&=&e^{zL(-1)}Y(b,-z)\psi_{V} \phi_{V}a\nonumber\\
&=&e^{zL(-1)}\psi_{V} \phi_{V}Y(\Delta(-h,-z)^{2}b,-z)a\nonumber\\
&=&e^{zL(-1)}\psi_{V} \phi_{V}e^{-zL(-1)}Y(a,z)\Delta(-h,-z)^{2}b\nonumber\\
&=&E^{-}(-h,z)^{2}\psi_{V} \phi_{V}Y(a,z)\Delta(-h,-z)^{2}b.
\end{eqnarray}
Similarly we obtain
\begin{eqnarray}
\bar{Y}(u,z)a&=&e^{zL(-1)}Y(a,-z)u\nonumber\\
&=&e^{zL(-1)}\phi_{V}Y(\Delta(-h,-z)a,-z)\phi_{V}^{-1} (u)\nonumber\\
&=&e^{zL(-1)}\phi_{V}e^{-zL(-1)}Y(\phi_{V}^{-1} (u),z)\Delta(-h,-z)a\nonumber\\
&=&e^{zL(-1)}\phi_{V}e^{-zL(-1)}\phi_{V}^{-1} Y(\Delta(h,z)\phi_{V}^{-1}
(u),z)\phi_{V}
\Delta(-h,-z)a\nonumber\\
&=&E^{-}(-h,z)Y(\Delta(h,z)\phi_{V}^{-1} (u),z)\phi_{V}\Delta(-h,-z)a.
\end{eqnarray}
For any $u,v\in \tilde{V}$, we get
\begin{eqnarray}
\bar{Y}(u,z)v&=&\phi_{V}^{-1} \bar{Y}(\Delta(h,z)u,z)\psi_{V} (v)\nonumber\\
&=&\phi_{V}^{-1} e^{zL(-1)}Y(\psi_{V}(v),-z)\Delta(h,z)u\nonumber\\
&=&\phi_{V}^{-1}
e^{zL(-1)}\phi_{V}Y(\Delta(-h,-z)\psi_{V}(v),-z)\phi_{V}^{-1}\Delta(h,z)u
\nonumber\\
&=&\phi_{V}^{-1} e^{zL(-1)}\phi_{V}e^{-zL(-1)}Y(\phi_{V}^{-1}\Delta(h,z)u,z)
\Delta(-h,-z)\psi_{V}(v)\nonumber\\
&=&E^{-}(-h,z)Y(\phi_{V}^{-1}\Delta(h,z)u,z)\Delta(-h,-z)\psi_{V}(v)\nonumber\\
&=&z^{-\gamma}E^{-}(-h,z)Y(\Delta(h,z)\phi_{V}^{-1}
(u),z)\Delta(-h,-z)\psi_{V}(v).
\end{eqnarray}
This proves the lemma. $\;\;\;\;\Box$

{\bf Lemma 3.6.} {\it The defined intertwining operator $\bar{Y}(\cdot,z)$ of
type
$\left(\begin{array}{c}V\\ \tilde{V} \tilde{V}\end{array}\right)$
satisfies the following condition:}
\begin{eqnarray}\label{e3.34}
\bar{Y}(u,z)v=(-1)^{\gamma}e^{zL(-1)}\bar{Y}(v,-z)u\;\;\;\mbox{{\it for any
}}u,v\in \tilde{V}.
\end{eqnarray}

{\bf Proof.} For any $u,v\in \tilde{V}$, using Lemma 3.5 we obtain
\begin{eqnarray}
& &\ \ \ e^{zL(-1)}\bar{Y}(v,-z)u\nonumber\\
& &=e^{zL(-1)}(-z)^{-\gamma}E^{-}(-h,-z)Y(\Delta(h,-z)\phi_{V}^{-1} (v),-z)
\Delta(-h,z)\psi_{V}(u)\nonumber\\
& &=e^{zL(-1)}(-z)^{-\gamma}E^{-}(-h,-z)Y(\Delta(h,-z)\phi_{V}^{-1} (v),-z)
\Delta(-h,z)\psi_{V}\phi_{V}\phi_{V}^{-1} (u)\nonumber\\
& &=(-z)^{-\gamma}z^{2\gamma}e^{zL(-1)}E^{-}(-h,-z)Y(\Delta(h,-z)
\phi_{V}^{-1} (v),-z)\psi_{V}\phi_{V}\Delta(-h,z)\phi_{V}^{-1} (u)\nonumber\\
&
&=(-z)^{\gamma}e^{zL(-1)}E^{-}(-h,-z)\psi_{V}\phi_{V}Y(\Delta(-h,-z)\phi_{V}^{-1}
(v),-z)\Delta(-h,z)\phi_{V}^{-1} (u)\nonumber\\
& &=(-z)^{\gamma}e^{zL(-1)}E^{-}(-h,-z)\psi_{V}\phi_{V}\nonumber\\
& &\ \ \ \cdot e^{-zL(-1)}Y(\Delta(-h,z)\phi_{V}^{-1}
(u),z)\Delta(-h,-z)\phi_{V}^{-1} (v)\nonumber\\
&
&=(-z)^{\gamma}(-z)^{-2\gamma}e^{zL(-1)}E^{-}(-h,-z)\psi_{V}\phi_{V}e^{-zL(-1)}\cdot \nonumber\\
& &\ \ \ \cdot Y(\Delta(-h,z)\phi_{V}^{-1}
(u),z)\phi_{V}^{-1}\psi_{V}^{-1}\Delta(-h,-z)\psi_{V} (v)\nonumber\\
& &=(-z)^{-\gamma}e^{zL(-1)}E^{-}(-h,-z)\psi_{V}\phi_{V}e^{-zL(-1)}
\phi_{V}^{-1}\psi_{V}^{-1}\nonumber\\
& &\ \ \ \cdot Y(\Delta(h,z)\phi_{V}^{-1} (u),z)\Delta(-h,-z)\psi_{V} (v)
\end{eqnarray}
and
\begin{eqnarray}
(-1)^{-\gamma}\bar{Y}(u,z)v
=(-z)^{-\gamma}E^{-}(-h,z)Y(\Delta(h,z)\phi_{V}^{-1}
(u),z)\Delta(-h,-z)\psi_{V} (v).
\end{eqnarray}
Using Lemma 3.4 we get
\begin{eqnarray}
&
&e^{zL(-1)}E^{-}(-h,-z)\psi_{V}\phi_{V}e^{-zL(-1)}\phi_{V}^{-1}\psi_{V}^{-1}\nonumber\\
&=&e^{zL(-1)}e^{-zL(-1)}\psi_{V} e^{zL(-1)}\phi_{V}e^{-zL(-1)}
\phi_{V}^{-1}\psi_{V}^{-1}\nonumber\\
&=&\psi_{V} e^{zL(-1)}\phi_{V}e^{-zL(-1)}
\phi_{V}^{-1}\psi_{V}^{-1}\nonumber\\
&=&\psi_{V} E^{-}(-h,z)\psi_{V}^{-1}\nonumber\\
&=&E^{-}(-h,z).
\end{eqnarray}
Then (\ref{e3.34}) follows immediately.$\;\;\;\;\Box$

{\bf Lemma 3.7.} {\it For any $a\in V$, the following identity hold}
\begin{eqnarray}
& &\ \ \ Y(E^{-}(h,z_{1})a,z_{2})\nonumber\\
& &=E^{-}(h,z_{1}+z_{2})E^{-}(-h,z_{2})Y(a,z_{2})z_{2}^{-h(0)}E^{+}(h,z_{2})
(z_{2}+z_{1})^{h(0)}E^{+}(-h,z_{2}+z_{1}).\nonumber\\
& &\mbox{}
\end{eqnarray}

{\bf Proof.} It is equivalent to prove that
\begin{eqnarray}
& &\ \ \ E^{-}(-h,z_{1}+z_{2})Y(E^{-}(h,z_{1})a,z_{2})(z_{2}+z_{1})^{-h(0)}
E^{+}(h,z_{2}+z_{1})\nonumber\\
& &=E^{-}(-h,z_{2})Y(a,z_{2})z_{2}^{-h(0)}E^{+}(h,z_{2}).
\end{eqnarray}
Since it is true when $z_{1}=0$, it is enough for us to prove that
partial derivatives for both sides with respect to the variable
$z_{1}$ are equal.
For any positive integer $k$ and for any $u\in V$, from the Jacobi identity
 we get
\begin{eqnarray}\label{e3.14}
& &Y(h(-k)u,z_{2})\nonumber\\
&=&\sum_{i=0}^{\infty}\left(\begin{array}{c}-k\\i\end{array}\right)\left(
(-z_{2})^{i}h(-k-i)Y(u,z_{2})-(-z_{2})^{-k-i}Y(u,z_{2})h(i)\right).
\end{eqnarray}
Using the binomial coefficient formulas:
$\left(\begin{array}{c}i-k\\i\end{array}\right)
=(-1)^{i}\left(\begin{array}{c}k-1\\i\end{array}\right)$ and
$\left(\begin{array}{c}-i\\k\end{array}\right)
=(-1)^{k+i-1}\left(\begin{array}{c}-k-1\\i-1\end{array}\right)$, then using
(\ref{e3.14})
for $u=E^{-}(h,z_{1})$ we obtain
\begin{eqnarray}
& &{\partial\over \partial z_{1}}Y(E^{-}(h,z_{1})a,z_{2})\nonumber\\
&=&Y\left(\sum_{k=1}^{\infty}z_{1}^{k-1}h(-k)E^{-}(h,z_{1})a,z_{2}\right)
\nonumber\\
&=&\sum_{k=1}^{\infty}z_{1}^{k-1}Y\left(h(-k)E^{-}(h,z_{1})a,z_{2}\right)
\nonumber\\
&=&\sum_{k=1}^{\infty}\sum_{i=0}^{\infty}\left(\begin{array}{c}-k\\i
\end{array}\right)z_{1}^{k-1}(-z_{2})^{i}h(-k-i)Y(u,z_{2})\nonumber\\
& &-\sum_{k=1}^{\infty}\sum_{i=0}^{\infty}\left(\begin{array}{c}-k\\i
\end{array}\right)z_{1}^{k-1}(-z_{2})^{-k-i}Y(u,z_{2})h(i)\nonumber\\
&=&\sum_{k=1}^{\infty}\sum_{i=0}^{\infty}\left(\begin{array}{c}i-k\\i
\end{array}\right)z_{1}^{k-i-1}(-z_{2})^{i}h(-k)Y(u,z_{2})\nonumber\\
& &-\sum_{k=1}^{\infty}z_{1}^{k-1}(-z_{2})^{-k}Y(u,z_{2})h(0)\nonumber\\
& &-\sum_{k=1}^{\infty}\sum_{i=1}^{\infty}\left(\begin{array}{c}-k\\i
\end{array}\right)z_{1}^{k-1}(-z_{2})^{-k-i}Y(u,z_{2})h(i)\nonumber\\
&=&\sum_{k=1}^{\infty}\sum_{i=0}^{\infty}(-1)^{i}\left(\begin{array}{c}k-1\\i
\end{array}\right)z_{1}^{k-i-1}(-z_{2})^{i}h(-k)Y(u,z_{2})\nonumber\\
& &+(z_{2}+z_{1})^{-1}Y(u,z_{2})h(0)
-\sum_{k=1}^{\infty}\sum_{i=1}^{\infty}\left(\begin{array}{c}-i\\k
\end{array}\right)z_{1}^{i-1}(-z_{2})^{-i-k}Y(u,z_{2})h(k)\nonumber\\
&=&\sum_{k=1}^{\infty}(z_{1}+z_{2})^{k-1}h(-k)Y(u,z_{2})
+(z_{2}+z_{1})^{-1}Y(u,z_{2})h(0)\nonumber\\
& &-\sum_{k=1}^{\infty}\sum_{i=1}^{\infty}(-1)^{k+i-1}
\left(\begin{array}{c}-k-1\\i-1
\end{array}\right)z_{1}^{i-1}(-z_{2})^{-i-k}Y(u,z_{2})h(k)\nonumber\\
&=&\sum_{k=1}^{\infty}(z_{1}+z_{2})^{k-1}h(-k)Y(u,z_{2})
+(z_{2}+z_{1})^{-1}Y(u,z_{2})h(0)\nonumber\\
& &+\sum_{k=1}^{\infty}(z_{2}+z_{1})^{-k-1}Y(u,z_{2})h(k).
\end{eqnarray}
Therefore
\begin{eqnarray}
{\partial\over \partial z_{1}}
E^{-}(-h,z_{1}+z_{2})Y(E^{-}(h,z_{1})a,z_{2})(z_{2}+z_{1})^{-h(0)}
E^{+}(h,z_{2}+z_{1})=0.
\end{eqnarray}
Then the proof is complete. $\;\;\;\;\Box$

{\bf Lemma 3.8.} {\it  For $a\in V$, we have}
\begin{eqnarray}
E^{-}(h,z_{1})Y(a,z_{2})E^{-}(-h,z_{1})
=Y(\Delta(-h,z_{2}-z_{1})\Delta(h,z_{2})a,z_{2}).
\end{eqnarray}

{\bf Proof.} For any $a\in V$, we have:
\begin{eqnarray}
[h(-k),
Y(a,z)]=\sum_{i=0}^{\infty}\left(\begin{array}{c}-k\\i\end{array}\right)z^{-k-i}
Y(h(i)a,z).
\end{eqnarray}
Noticing that $\displaystyle{{1\over
k}\left(\begin{array}{c}-k\\i\end{array}\right)=(-1)^{i+k}
\left(\begin{array}{c}-i\\k\end{array}\right)}$ for any nonnegative integers
$k$ and $i$, we
have:
\begin{eqnarray}
& &\left[\sum_{k=1}^{\infty}\frac{h(-k)}{k}z_{1}^{k}, Y(a,z_{2})\right]
\nonumber\\
&=&\sum_{k=1}^{\infty}\sum_{i=0}^{\infty}{1\over k}
\left(\begin{array}{c}-k\\i\end{array}\right)
z_{1}^{k}z_{2}^{-k-i}Y(h(i)a,z_{2})\nonumber\\
&=&\sum_{k=1}^{\infty}{1\over k}z_{1}^{k}z_{2}^{-k}Y(h(0)a,z_{2})
+\sum_{k=1}^{\infty}\sum_{i=1}^{\infty}{1\over k}
\left(\begin{array}{c}-k\\i\end{array}\right)
z_{1}^{k}z_{2}^{-k-i}Y(h(i)a,z_{2})\nonumber\\
&=&-\log \left(1-{z_{1}\over z_{2}}\right)Y(h(0)a,z_{2})\nonumber\\
& &+\sum_{k=1}^{\infty}\sum_{i=1}^{\infty}(-1)^{i+k}{1\over i}
\left(\begin{array}{c}-i\\k\end{array}\right)
z_{1}^{k}z_{2}^{-k-i}Y(h(i)a,z_{2})\nonumber\\
&=&-\log \left(1-{z_{1}\over z_{2}}\right)Y(h(0)a,z_{2})\nonumber\\
& &+\sum_{i=1}^{\infty}{1\over i}\left((-z_{2}+z_{1})^{-i}-(-z_{2})^{-i}\right)
Y(h(i)a,z_{2}).
\end{eqnarray}
Then
\begin{eqnarray}
& &E^{-}(h,z_{1})Y(a,z_{2})E^{-}(-h,z_{1})\nonumber\\
&=&Y\left(\left(1-{z_{1}\over z_{2}}\right)^{-h(0)}E^{+}(h,-z_{2}+z_{1})
E^{+}(-h,-z_{2})a,z_{2}\right)\nonumber\\
&=&Y\left(\Delta(-h,z_{2}-z_{1})\Delta(h,z_{2})a,z_{2}\right).
\;\;\;\;\Box
\end{eqnarray}

Now we are ready to prove our main theorem.

{\bf Theorem 3.9.} {\it Let $V$ be a vertex operator algebra and let $h\in V$
satisfying the conditions (\ref{e2.15}) with an
integer $\gamma$. Then
a) if $\gamma$ is even, $V\oplus \tilde{V}$ is a vertex operator
algebra. b) if $\gamma$ is odd, $V\oplus \tilde{V}$
is a vertex operator superalgebra.}

{\bf Proof.} From [FHL], we only need to prove the Jacobi identity
for three module elements. For any $u,v,w\in \tilde{V}$, since
$\bar{Y}(v,z)w\in V((z))$,
using Lemmas 3.5 and 3.8 we obtain
\begin{eqnarray}
& &\ \ \ \bar{Y}(u,z_{1})\bar{Y}(v,z_{2})w\nonumber\\
& &=E^{-}(-h,z_{1})Y(\Delta(h,z_{1})\phi_{V}^{-1}
(u),z_{1})\phi_{V}\Delta(-h,-z_{1})
\bar{Y}(v,z_{2})w\nonumber\\
& &=E^{-}(-h,z_{1})Y(\Delta(h,z_{1})\phi_{V}^{-1}
(u),z_{1})\phi_{V}\Delta(-h,-z_{1})
\cdot\nonumber\\
& &\ \ \ \cdot
z_{2}^{-\gamma}E^{-}(-h,z_{2})Y(\Delta(h,z_{2})\phi_{V}^{-1}v,z_{2})\Delta(-h,-z_{2})\psi_{V}(w)
\nonumber\\
& &=E^{-}(-h,z_{1})Y(\Delta(h,z_{1})\phi_{V}^{-1} (u),z_{1})E^{-}(-h,z_{2})
\left(1-{z_{2}\over z_{1}}\right)^{\gamma}z_{2}^{-\gamma}
\cdot\nonumber\\
& &\ \ \ \cdot \phi_{V}\Delta(-h,-z_{1})
Y(\Delta(h,z_{2})\phi_{V}^{-1}v,z_{2})\Delta(-h,-z_{2})\psi_{V}(w)\nonumber\\
& &=E^{-}(-h,z_{1})E^{-}(-h,z_{2})
Y(\Delta(-h,z_{1}-z_{2})\Delta(h,z_{1})^{2}\phi_{V}^{-1} (u),z_{1})
\left(1-{z_{2}\over z_{1}}\right)^{\gamma}z_{2}^{-\gamma}\cdot\nonumber\\
& &\ \ \ \cdot \phi_{V}\Delta(-h,-z_{1})
Y(\Delta(h,z_{2})\phi_{V}^{-1}v,z_{2})\Delta(-h,-z_{2})\psi_{V}(w)\nonumber\\
& &=E^{-}(-h,z_{1})E^{-}(-h,z_{2})
Y(\Delta(-h,z_{1}-z_{2})\Delta(h,z_{1})^{2}\phi_{V}^{-1} (u),z_{1})
\left(1-{z_{2}\over z_{1}}\right)^{\gamma}z_{2}^{-\gamma}\cdot\nonumber\\
& &\ \ \ \cdot \phi_{V}
Y(\Delta(-h,-z_{1}+z_{2})\Delta(h,z_{2})\phi_{V}^{-1}v,z_{2})
\Delta(-h,-z_{1})\Delta(-h,-z_{2})\psi_{V}(w)\nonumber\\
& &=E^{-}(-h,z_{1})E^{-}(-h,z_{2})
Y(\Delta(-h,z_{1}-z_{2})\Delta(h,z_{1})^{2}\phi_{V}^{-1}u,z_{1})
\left(1-{z_{2}\over z_{1}}\right)^{\gamma}z_{2}^{-\gamma}(-z_{1})^{\gamma}
\cdot\nonumber\\
& &\ \ \ \cdot (-z_{2})^{\gamma}
\phi_{V}Y(\Delta(-h,-z_{1}+z_{2})\Delta(h,z_{2})\phi_{V}^{-1}v,z_{2})
\psi_{V}\Delta(-h,-z_{1})\Delta(-h,-z_{2})w\nonumber\\
& &=\phi_{V}\psi_{V} E^{-}(-h,z_{1})E^{-}(-h,z_{2})
Y(\Delta(-h,z_{1}-z_{2})\phi_{V}^{-1} (u),z_{1})
(z_{1}-z_{2})^{\gamma}\cdot\nonumber\\
& &\ \ \ \cdot Y(\Delta(-h,-z_{1}+z_{2})\phi_{V}^{-1}v,z_{2})
\Delta(-h,-z_{1})\Delta(-h,-z_{2})w.
\end{eqnarray}
Symmetrically, we have:
\begin{eqnarray}
& &\bar{Y}(v,z_{2})\bar{Y}(u,z_{1})w\nonumber\\
&=&\phi_{V}\psi_{V} E^{-}(-h,z_{1})E^{-}(-h,z_{2})
Y(\Delta(-h,z_{2}-z_{1})\phi_{V}^{-1}v,z_{2})
(z_{2}-z_{1})^{\gamma}\cdot\nonumber\\
& &\cdot Y(\Delta(-h,-z_{2}+z_{1})\phi_{V}^{-1}u,z_{1})
\Delta(-h,-z_{1})\Delta(-h,-z_{2})w.
\end{eqnarray}
Then
\begin{eqnarray}
& &z_{0}^{-1}\delta\left(\frac{z_{1}-z_{2}}{z_{0}}\right)
\bar{Y}(u,z_{1})\bar{Y}(v,z_{2})w
-(-1)^{\gamma}z_{0}^{-1}\delta\left(\frac{z_{2}-z_{1}}{-z_{0}}\right)
\bar{Y}(v,z_{2})\bar{Y}(u,z_{1})w
\nonumber\\
&=&z_{0}^{-1}\delta\left(\frac{z_{1}-z_{2}}{z_{0}}\right)
\phi_{V}\psi_{V} E^{-}(-h,z_{1})E^{-}(-h,z_{2})z_{0}^{\gamma}\cdot \nonumber\\
& &\cdot
Y(\Delta(-h,z_{0})\phi_{V}^{-1}u,z_{1})Y(\Delta(-h,-z_{0})\phi_{V}^{-1}v,z_{2})C
\nonumber\\
& &-z_{0}^{-1}\delta\left(\frac{z_{2}-z_{1}}{-z_{0}}\right)
\phi_{V}\psi_{V} E^{-}(-h,z_{1})E^{-}(-h,z_{2})
z_{0}^{\gamma}\cdot \nonumber\\
& &\cdot Y(\Delta(-h,-z_{0})\phi_{V}^{-1}
(v),z_{2})Y(\Delta(-h,z_{0})\phi_{V}^{-1}u,z_{1})C
\nonumber\\
&=&z_{2}^{-1}\delta\left(\frac{z_{1}-z_{0}}{z_{2}}\right)
\phi_{V}\psi_{V} E^{-}(-h,z_{1})E^{-}(-h,z_{2})z_{0}
^{\gamma}Y(Y(A,z_{0})B,z_{2})C,
\end{eqnarray}
where
\begin{eqnarray}
A=\Delta(-h,z_{0})\phi_{V}^{-1} (u),\;
B=\Delta(-h,-z_{0})\phi_{V}^{-1} (v),\;
C=\Delta(-h,-z_{1})\Delta(-h,-z_{2})w.
\end{eqnarray}
On the other hand, using Lemmas 3.5 and 3.7 we obtain
\begin{eqnarray}
&
&z_{2}^{-1}\delta\left(\frac{z_{1}-z_{0}}{z_{2}}\right)\bar{Y}(\bar{Y}(u,z_{0})v,z_{2})w
\nonumber\\
&=&z_{2}^{-1}\delta\left(\frac{z_{1}-z_{0}}{z_{2}}\right)
Y\left(E^{-}(-h,z_{0})Y(\phi_{V}^{-1}\Delta(h,z_{0})u,z_{0})\Delta(-h,-z_{0})\psi_{V}(v),
z_{2}\right)w\nonumber\\
&=&z_{2}^{-1}\delta\left(\frac{z_{1}-z_{0}}{z_{2}}\right)
E^{-}(-h,z_{0}+z_{2})E^{-}(h,z_{2})\nonumber\\
& &\cdot Y\left(Y(\phi_{V}^{-1}\Delta(h,z_{0})u,z_{0})
\Delta(-h,-z_{0})\psi_{V}(v),z_{2}\right)\cdot\nonumber\\
& &\cdot z_{2}^{h(0)}E^{+}(-h,z_{2})(z_{2}+z_{0})^{-h(0)}
E^{+}(h,z_{2}+z_{0})w\nonumber\\
&=&z_{2}^{-1}\delta\left(\frac{z_{1}-z_{0}}{z_{2}}\right)
E^{-}(-h,z_{0}+z_{2})E^{-}(h,z_{2})\cdot\nonumber\\
& &\cdot
Y\left(Y(\phi_{V}^{-1}\Delta(h,z_{0})u,z_{0})\Delta(-h,-z_{0})\psi_{V}(v),
z_{2}\right)\Delta(h,-z_{2})\Delta(-h,-z_{2}-z_{0})w\nonumber\\
&=&z_{2}^{-1}\delta\left(\frac{z_{1}-z_{0}}{z_{2}}\right)
E^{-}(-h,z_{1})E^{-}(h,z_{2})z_{0}^{-\gamma}\cdot\nonumber\\
& &\cdot Y\left(Y(\Delta(h,z_{0})\phi_{V}^{-1}
(u),z_{0})\Delta(-h,-z_{0})\psi_{V}(v),
z_{2}\right)\Delta(h,-z_{2})\Delta(-h,-z_{1})w\nonumber\\
&=&z_{2}^{-1}\delta\left(\frac{z_{1}-z_{0}}{z_{2}}\right)
E^{-}(-h,z_{1})E^{-}(h,z_{2})z_{0}^{-\gamma}(-z_{0})^{2\gamma}\cdot\nonumber\\
& &\cdot Y\left(Y(\Delta(h,z_{0})\phi_{V}^{-1} (u),z_{0})\psi_{V}
\phi_{V}\Delta(-h,-z_{0})\phi_{V}^{-1} (v),
z_{2}\right)\Delta(h,-z_{2})\Delta(-h,-z_{1})w\nonumber\\
&=&z_{2}^{-1}\delta\left(\frac{z_{1}-z_{0}}{z_{2}}\right)
E^{-}(-h,z_{1})E^{-}(h,z_{2})z_{0}^{\gamma}\cdot\nonumber\\
& &\cdot Y\left(\psi_{V}\phi_{V} Y(\Delta(-h,z_{0})\phi_{V}^{-1} (u),z_{0})
\Delta(-h,-z_{0})\phi_{V}^{-1} (v),
z_{2}\right)\Delta(h,-z_{2})\Delta(-h,-z_{1})w\nonumber\\
&=&z_{2}^{-1}\delta\left(\frac{z_{1}-z_{0}}{z_{2}}\right)
\psi_{V}\phi_{V}E^{-}(-h,z_{1})E^{-}(-h,z_{2})z_{0}^{\gamma}\cdot\nonumber\\
& &\cdot Y\left( Y(\Delta(-h,z_{0})\phi_{V}^{-1} (u),z_{0})
\Delta(-h,-z_{0})\phi_{V}^{-1} (v),
z_{2}\right)\Delta(-h,-z_{2})\Delta(-h,-z_{1})w.
\end{eqnarray}
Therefore, the Jacobi identity holds. Then the proof is complete.$\;\;\;\;\Box$

{\bf Remark 3.10.} Under the assumption of Theorem 3.9, suppose that $V$ is a
simple vertex
operator algebra. Then $\bar{V}=V\oplus \tilde{V}$ is a completely reducible
$V$-module.
Let $P_{0}$ and $P_{1}$ be the projection maps from
$\bar{V}$ onto $V$ and $\tilde{V}$, respectively. Then $P_{0}$ and $P_{1}$ are
$V$-homomorphisms. Let $W$ be any irreducible $V$-submodule of $\bar{V}$.
Then $P_{0}$ and
$P_{1}$ are  $V$-homomorphisms from $W$ to $V$ and $\tilde{V}$, respectively.
Thus, if $\tilde{V}$ as a $V$-module is not isomorphic to
 $V$, either $P_{0}(W)=0$ or $P_{1}(W)=0$, so that $W=V$ or $W=\bar{V}$.
Therefore, if
$\tilde{V}$ as a $V$-module is not isomorphic to
 $V$, $\bar{V}$ is a simple vertex operator superalgebra.

Next, continuing with Theorem 3.9  we study modules for the vertex operator
(super)algebra $\bar{V}$. We define a linear endomorphism $\sigma$ of
$\bar{V}$ by: $\sigma|_{V}=1,\; \sigma|_{\tilde{V}}=-1$. Then it is
 clear that $\sigma$ is an order-two automorphism of the vertex operator
algebra $\bar{V}$. Recall that $\pi_{0}$ is a $V$-isomorphism from
$V$ onto $\tilde{\tilde{V}}$. Let $W$ be any $V$-module. Then we shall show
that $\tilde{\tilde{W}}$ is isomorphic to $W$ as a $V$-module. For simplicity
we shall
prove this only for an irreducible $V$-module $W$.

{\bf Proposition 3.11.} {\it Let $W$ be an irreducible $V$-module. Then the
linear map}
\begin{eqnarray}
I(\cdot,z):& & V\rightarrow {\rm Hom}(W,\tilde{\tilde{W}})\{z\};\nonumber\\
& &I(a,z)w=e^{zL(-1)}\psi_{\tilde{W}}^{-1}\psi_{W}^{-1}e^{-zL(-1)}
Y_{W}(\psi_{V}\psi_{\tilde{V}}\pi_{0}(a),z)\Delta(2h,-z)w.
\end{eqnarray}
{\it is a nonzero intertwining operator of type $\left(\begin{array}{c}
\tilde{\tilde{W}}\\V W\end{array}\right)$. Consequently,
 $\pi_{W}:=I({\bf 1},z)$ is a  $V$-isomorphism
from $W$ onto $\tilde{\tilde{W}}$.}

{\bf Proof.} We shall construct $I(\cdot,z)$ from the intertwining operator
$Y_{W}(\cdot,z)$
through the following sequence:
\begin{eqnarray}
\left(\begin{array}{c}W\\V W\end{array}\right)\rightarrow
\left(\begin{array}{c}W\\W V\end{array}\right)\rightarrow
\left(\begin{array}{c}\tilde{\tilde{W}}\\W \tilde{\tilde{V}}\end{array}\right)
\rightarrow
\left(\begin{array}{c}\tilde{\tilde{W}}\\W V\end{array}\right)
\rightarrow
\left(\begin{array}{c}\tilde{\tilde{W}}\\V W\end{array}\right).
\end{eqnarray}
Starting from the intertwining operator $Y_{W}(\cdot,z)$ of type
$\left(\begin{array}{c}W\\V W\end{array}\right)$ we have an intertwining
operator $I_{1}(\cdot,z)$ defined as follows [FHL]:
\begin{eqnarray}
I_{1}(w,z)a=e^{zL(-1)}Y_{W}(a,-z)w\;\;\;\mbox{for any }a\in V,w\in W.
\end{eqnarray}
Since $Y_{W}(\cdot,z)$ is not zero, $I_{1}(\cdot,z)$ is not zero either.
By Lemma 2.1 we have an intertwining operator $I_{2}(\cdot,z)$ of type
$\left(\begin{array}{c}\tilde{\tilde{W}}\\W \tilde{\tilde{V}}
\end{array}\right)$ defined by
\begin{eqnarray}
I_{2}(w,z)u=\psi_{\tilde{W}}^{-1}\psi_{W}^{-1}I_{1}(\Delta(2h,z)w,z)\psi_{V}\psi_{\tilde{V}}u
\;\;\;\mbox{for any }u\in \tilde{\tilde{V}}, w\in W.
\end{eqnarray}
Using the isomorphism $\pi_{0}$ from $\tilde{\tilde{V}}$ onto $V$ we obtain
an intertwining operator $I_{3}(\cdot,z)$ of type
$\left(\begin{array}{c}\tilde{\tilde{W}}\\W V\end{array}\right)$ defined as
follows:
\begin{eqnarray}
I_{3}(w,z)a=I_{2}(w,z)\pi_{0}a\;\;\;\mbox{for any }w\in W,a\in V.
\end{eqnarray}
Finally we obtain an intertwining operator $I(\cdot,z)$ of type
$\left(\begin{array}{c}\tilde{\tilde{W}}\\V W\end{array}\right)$ defined by
\begin{eqnarray}
I(a,z)w=e^{zL(-1)}I_{3}(w,-z)a\;\;\;\mbox{for any }a\in V,w\in W.
\end{eqnarray}
Composing all intertwining operators together we obtain
\begin{eqnarray}
I(a,z)w&=&e^{zL(-1)}I_{3}(w,-z)a\nonumber\\
&=&e^{zL(-1)}I_{2}(w,-z)\pi_{0}(a)\nonumber\\
&=&e^{zL(-1)}\psi_{\tilde{W}}^{-1}\psi_{W}^{-1}I_{1}(\Delta(2h,-z)w,-z)
\psi_{V}\psi_{\tilde{V}}\pi_{0}(a)\nonumber\\
&=&e^{zL(-1)}\psi_{\tilde{W}}^{-1}\psi_{W}^{-1}e^{-zL(-1)}
Y_{W}(\psi_{V}\psi_{\tilde{V}}\pi_{0}(a),z)\Delta(2h,-z)w.
\end{eqnarray}
Because
\begin{eqnarray}
{d\over dz}I({\bf 1},z)=I(L(-1){\bf 1},z)=0,
\end{eqnarray}
$I({\bf 1},z)$ is a constant linear map from $W$ to $\tilde{\tilde{W}}$.
Since $[Y(a,z_{1}), I({\bf 1},z)]=0$ for any $a\in V$, $I({\bf 1},z)$ is a
$V$-homomorphism from $W$ to $\tilde{\tilde{W}}$. Since $W$ and
$\tilde{\tilde{W}}$ are
irreducible $V$-modules, $I({\bf 1},z)$ is a $V$-isomorphism.
Then the proof is complete. $\;\;\;\;\Box$

Let $W$ be an irreducible $V$-module (with finite-dimensional homogeneous
subspaces). Since $h(0)$ preserves each homogeneous subspace of $M$ and $h(0)$
is
semisimple on $V$, $h(0)$ acts semisimply on $W$.
{\em Suppose that $h(0)$ has only rational  eigenvalues on $W$}.
Set
$\phi_{W}=\psi_{\tilde{W}}\pi_{W}$. That is,
$\phi_{W}$ is a linear isomorphism from $W$ onto $\tilde{W}$
such that
\begin{eqnarray}
\phi_{W}(Y(a,z)w)=Y(\Delta(h,z)a,z)\phi_{W}(w)\;\;\;\mbox{ for}a\in V, w\in W.
\end{eqnarray}

Next we shall make $\bar{W}:=W\oplus \tilde{W}$ a $\tilde{V}$-module
or $\sigma$-twisted $\tilde{V}$-module. In order to do this we need
intertwining
operators of types $\left(\begin{array}{c}\tilde{W}\\ \tilde{V}
W\end{array}\right)$
and $\left(\begin{array}{c}W\\ \tilde{V} \tilde{W}\end{array}\right)$.

Following [DL], for any $\alpha\in {\bf Q}$ we define
$(-z)^{\alpha}=e^{\alpha\pi i}$.
For any $a\in V, w\in W$,  define $F(u,z)a=e^{zL(-1)}Y_{W}(a,-z)u$. Then
$F(\cdot,z)$ is an intertwining operator of type
$\left(\begin{array}{c}W\\ W V\end{array}\right)$. For any $u\in \tilde{V},
w\in W$, define
$F_{1}(w,z)u=\psi_{W}^{-1}F(\Delta(h,z)w,z)\psi_{V}(u)$. Then from Lemma 2.1,
$F_{1}(\cdot,z)$ is an intertwining operator of type
$\left(\begin{array}{c}\tilde{W}\\ W \tilde{V}\end{array}\right)$.
For any $u\in \tilde{V}, w\in W$ we define $Y(u,z)w=e^{zL(-1)}F_{1}(w,-z)u$.
Then
$Y(\cdot,z)$ is an intertwining operator of type
$\left(\begin{array}{c}\tilde{W}\\ \tilde{V} W\end{array}\right)$.
That is,
\begin{eqnarray}
Y(u,z)w
&=&e^{zL(-1)}F_{1}(w,-z)u\nonumber\\
&=&e^{zL(-1)}\psi_{W}^{-1}F(\Delta(h,-z)w,-z)\psi_{V}(u)
\nonumber\\
&=&e^{zL(-1)}\psi_{W}^{-1}e^{-zL(-1)}Y_{W}(\psi_{V}(u),z)\Delta(h,-z)w
\nonumber\\
&=&E^{-}(h,z)\psi_{W}^{-1}Y(\psi_{V}(u),z)\Delta(h,-z)w.
\end{eqnarray}
Then we obtain an intertwining operator of type
$\left(\begin{array}{c}\tilde{W}\\ \tilde{V} W\end{array}\right)$.

{}From Lemma 2.1 we obtain an intertwining operator $\tilde{Y}(\cdot,z)$ of
type
$\left(\begin{array}{c}\tilde{\tilde{W}}\\ \tilde{V}
\tilde{W}\end{array}\right)$.
Then $Y(\cdot,z):=\pi_{W}^{-1}\tilde{Y}(\cdot,z)$ is an intertwining operator
of type
$\left(\begin{array}{c}W\\ \tilde{V} \tilde{W}\end{array}\right)$.
For any
$u\in \tilde{V}, \tilde{w}\in \tilde{W}$, we have:
\begin{eqnarray}
Y(u,z)\tilde{w}
&=&\pi_{W}^{-1}\tilde{Y}(u,z)\tilde{w}\nonumber\\
&=&\pi_{W}^{-1}\psi_{\tilde{W}}^{-1}Y(\Delta(h,z)u,z)\psi_{W}\tilde{w}\nonumber\\
&=&E^{-}(h,z)\pi_{W}^{-1}\psi_{\tilde{W}}^{-1}\psi_{W}^{-1}
Y(\psi_{V}\Delta(h,z)(u),z)\Delta(h,-z)\psi_{W}\tilde{w}.
\end{eqnarray}
Then we obtain an intertwining operator of type
$\left(\begin{array}{c}W\\ \tilde{V} \tilde{W}\end{array}\right)$.

Similar to Lemma 3.6 we have:

{\bf Lemma 3.12.} {\it For $a\in V, u\in \tilde{V},w\in W, \tilde{w}\in
\tilde{W}$, we have}
\begin{eqnarray}
& &Y(\psi_{V}\phi_{V}
a,z)u=E^{-}(-2h,z)\psi_{W}\psi_{\tilde{W}}\pi_{W}Y(a,z)\Delta(-2h,-z)u,\\
& &Y(u,z)w=
E^{-}(-h,z)\phi_{W}Y(\phi_{V}^{-1}u,z)\Delta(-h,-z)w,\\
& &Y(u,z)\tilde{w}
=(-1)^{\gamma}E^{-}(-h,z)\psi_{W}Y(\phi_{V}^{-1}u,z)
\Delta(-h,-z)\tilde{w}.
\end{eqnarray}

{\bf Proof.} Recall from Proposition 3.11 that $I(\cdot,z)$ is an intertwining
operator
of type $\left(\begin{array}{c}
\tilde{\tilde{W}}\\V W\end{array}\right)$. Then $\pi_{W}^{-1}I(\cdot,z)$ is an
intertwining operator
of type $\left(\begin{array}{c}W\\V W\end{array}\right)$. Since
$I\left(\begin{array}{c}W\\V W\end{array}\right)$ is one-dimensional and
$\pi_{W}^{-1}I({\bf 1},z)=id_{W}=Y_{W}({\bf 1},z)$, we get
$\pi_{W}^{-1}I(\cdot,z)=Y_{W}(\cdot,z)$.
Thus
\begin{eqnarray}
Y_{W}(a,z)w=\pi_{W}^{-1}e^{zL(-1)}\psi_{\tilde{W}}^{-1}\psi_{W}^{-1}e^{-zL(-1)}
Y_{W}(\psi_{V}\psi_{\tilde{V}}\pi_{0}(a),z)\Delta(2h,-z)w.
\end{eqnarray}
Therefore (by Lemma 3.4)
\begin{eqnarray}
& &Y_{W}(\psi_{V}\psi_{\tilde{V}}\pi_{0}(a),z)w\nonumber\\
&=&e^{zL(-1)}\psi_{W}\psi_{\tilde{W}}e^{-zL(-1)}
\pi_{W}Y_{W}(a,z)\Delta(-2h,-z)w\nonumber\\
&=&E^{-}(-2h,z)\psi_{W}\psi_{\tilde{W}}\pi_{W}Y_{W}(a,z)\Delta(-2h,-z)w.
\end{eqnarray}
For any $u\in \tilde{V}, w\in W$, using the first identity we obtain
\begin{eqnarray}
Y(u,z)w
&=&E^{-}(h,z)\psi_{W}^{-1}Y(\psi_{V}\phi_{V}\phi_{V}^{-1}u,z)\Delta(h,-z)w\nonumber\\
&=&E^{-}(-h,z)\psi_{\tilde{W}}\pi_{W}Y(\phi_{V}^{-1}u,z)\Delta(-h,-z)w.
\end{eqnarray}
Similarly, for any
$u\in \tilde{V}, \tilde{w}\in \tilde{W}$, we obtain
\begin{eqnarray}
Y(u,z)\tilde{w}
&=&E^{-}(-h,z)Y(\phi_{V}^{-1}\Delta(h,z)u,z)\Delta(-h,-z)\psi_{W}\tilde{w}\nonumber\\
&=&z^{-\gamma}E^{-}(-h,z)Y(\Delta(h,z)\phi_{V}^{-1}u,z)\Delta(-h,-z)\psi_{W}\tilde{w}
\nonumber\\
&=&(-1)^{\gamma}E^{-}(-h,z)\psi_{W}Y(\phi_{V}^{-1}u,z)\Delta(-h,-z)\tilde{w}.
\end{eqnarray}
Then the proof is complete.$\;\;\;\;\Box$

Recall ([D2], [FFR], [FLM])  that a $\sigma$-twisted $\bar{V}$-module is a
${\bf C}$-graded
vector space $M=\oplus_{h\in {\bf C}}M_{h}$ satisfying all conditions for a
module except that
$Y_{M}(a,z)\!\in\! {\rm End}M[[z^{{1\over 2}}, z^{-{1\over 2}}]]$ for $a\in
\bar{V}$ and
that
the Jacobi identity is replaced by the following twisted Jacobi identity:
 \begin{equation}\label{jacobi}
\begin{array}{c}
\displaystyle{z^{-1}_0\delta\left(\frac{z_1-z_2}{z_0}\right)
Y_M(u,z_1)Y_M(v,z_2)-z^{-1}_0\delta\left(\frac{z_2-z_1}{-z_0}\right)
Y_M(v,z_2)Y_M(u,z_1)}\\
\displaystyle{=z_2^{-1}\left(\frac{z_1-z_0}{z_2}\right)^{-r/2}
\delta\left(\frac{z_1-z_0}{z_2}\right)
Y_M(Y(u,z_0)v,z_2)}
\end{array}
\end{equation}
for $u\in V^r$, where $r=0,1$ and $V^0=V, V^1=\tilde{V}$.

{\bf Theorem 3.13. } {\it Let $W$ be an irreducible $V$-module such that
$h(0)$ has only half-integral eigenvalues on $W$. Then
$W\oplus \tilde{W}$ is either a $\bar{V}$-module or a
 $\sigma$-twisted $\bar{V}$-module.}

{\bf Proof.} It is similar to the proof of Theorem 3.9.
We only need to prove the Jacobi identity for $(u,v,w)$, where
$u,v\in \tilde{V}, w\in W\cup \tilde{W}$.

Case 1: $w\in W$. By Lemma 3.12 we obtain
\begin{eqnarray}
& &Y(u,z_{1})Y(v,z_{2})w\nonumber\\
&=&(-1)^{\gamma}E^{-}(-h,z_{1})\psi_{W}Y(\phi_{V}^{-1}u,
z_{1})\Delta(-h,-z_{1})\cdot\nonumber\\
& &\cdot
E^{-}(-h,z_{2})\phi_{W}Y(\phi_{V}^{-1}v,z_{2})\Delta(-h,-z_{2})w\nonumber\\
&=&(-1)^{\gamma}\left(1-{z_{2}\over z_{1}}\right)^{\gamma}E^{-}(-h,z_{1})
\psi_{W}Y(\phi_{V}^{-1}u,z_{1})\cdot \nonumber\\
& &\cdot E^{-}(-h,z_{2})\Delta(-h,-z_{1})\phi_{W}Y(\phi_{V}^{-1}v,z_{2})
\Delta(-h,-z_{2})w\nonumber\\
&=&(-1)^{\gamma}\left(1-{z_{2}\over z_{1}}\right)^{\gamma}
E^{-}(-h,z_{1})E^{-}(-h,z_{2})\psi_{W}
Y(\Delta(-h,z_{1}-z_{2})\Delta(h,z_{1})\phi_{V}^{-1}u,z_{1})\cdot \nonumber\\
& &\cdot \Delta(-h,-z_{1})\phi_{W}Y(\phi_{V}^{-1}v,z_{2})
\Delta(-h,-z_{2})w\nonumber\\
&=&(-1)^{\gamma}\left(1-{z_{2}\over z_{1}}\right)^{\gamma}
E^{-}(-h,z_{1})E^{-}(-h,z_{2})\psi_{W}
Y(\Delta(-h,z_{1}-z_{2})\Delta(h,z_{1})\phi_{V}^{-1}u,z_{1})\cdot \nonumber\\
& &\cdot (-z_{1})^{\gamma}\phi_{W}\Delta(-h,-z_{1})Y(\phi_{V}^{-1}v,z_{2})
\Delta(-h,-z_{2})w\nonumber\\
&=&(z_{1}-z_{2})^{\gamma}
E^{-}(-h,z_{1})E^{-}(-h,z_{2})\psi_{W}\phi_{W}
Y(\Delta(-h,z_{1}-z_{2})\phi_{V}^{-1}u,z_{1})\cdot \nonumber\\
& &\cdot Y(\Delta(-h,-z_{1}+z_{2})\phi_{V}^{-1}v,z_{2})
\Delta(-h,-z_{1})\Delta(-h,-z_{2})w.
\end{eqnarray}
Then
\begin{eqnarray}
& &z_{0}^{-1}\delta\left(\frac{z_{1}-z_{2}}{z_{0}}\right)
Y(u,z_{1})Y(v,z_{2})w\nonumber\\
&=&z_{0}^{-1}\delta\left(\frac{z_{1}-z_{2}}{z_{0}}\right)z_{0}^{\gamma}
E^{-}(-h,z_{1})E^{-}(-h,z_{2})\psi_{W}\phi_{W}
Y(\Delta(-h,z_{0})\phi_{V}^{-1} (u),z_{1})\cdot\nonumber\\
& &\cdot Y(\Delta(-h,-z_{0})\phi_{V}^{-1} (v),z_{2})
\Delta(-h,-z_{1})\Delta(-h,-z_{2})w.
\end{eqnarray}
Symmetrically we have:
\begin{eqnarray}
& &z_{0}^{-1}\delta\left(\frac{z_{2}-z_{1}}{-z_{0}}\right)
Y(v,z_{2})Y(u,z_{1})w\nonumber\\
&=&z_{0}^{-1}\delta\left(\frac{z_{1}-z_{2}}{z_{0}}\right)(-z_{0})^{\gamma}
E^{-}(-h,z_{1})E^{-}(-h,z_{2})\psi_{W}\phi_{W}
 Y(\Delta(-h,-z_{0})\phi_{V}^{-1} (v),z_{2})\cdot \nonumber\\
& &\cdot Y(\Delta(-h,z_{0})\phi_{V}^{-1} (u),z_{1})
\Delta(-h,-z_{1})\Delta(-h,-z_{2})w.
\end{eqnarray}
Therefore
\begin{eqnarray}
& &z_{0}^{-1}\delta\left(\frac{z_{1}-z_{2}}{z_{0}}\right)
Y(u,z_{1})Y(v,z_{2})w
-(-1)^{\gamma}z_{0}^{-1}\delta\left(\frac{z_{2}-z_{1}}{-z_{0}}\right)
Y(v,z_{2})Y(u,z_{1})w\nonumber\\
&=&z_{2}^{-1}\delta\left(\frac{z_{1}-z_{0}}{z_{2}}\right)z_{0}^{\gamma}
E^{-}(-h,z_{1})E^{-}(-h,z_{2})\psi_{W}\phi_{W}\cdot \nonumber\\
& &\cdot Y(Y(\Delta(-h,z_{0})\phi_{V}^{-1} (u),z_{0})
\Delta(-h,-z_{0})\phi_{V}^{-1} (v),z_{2})C
\end{eqnarray}
where $C=\Delta(-h,-z_{1})\Delta(-h,-z_{2})w$.
On the other hand, we have:
\begin{eqnarray}
& &z_{2}^{-1}\delta\left(\frac{z_{1}-z_{0}}{z_{2}}\right)
Y(Y(u,z_{0})v,z_{2})w\nonumber\\
&=&z_{2}^{-1}\delta\left(\frac{z_{1}-z_{0}}{z_{2}}\right)z_{0}^{-\gamma}
Y\left(E^{-}(-h,z_{0})Y(\Delta(h,z_{0})\phi_{V}^{-1}u,z_{0})\Delta(-h,-z_{0})\psi_{V}(v),
z_{2}\right)w\nonumber\\
&=&z_{2}^{-1}\delta\left(\frac{z_{1}-z_{0}}{z_{2}}\right)
E^{-}(-h,z_{1})E^{-}(h,z_{2})z_{0}^{-\gamma}\cdot\nonumber\\
& &\cdot Y\left(Y(\Delta(h,z_{0})\phi_{V}^{-1} (u),z_{0})
\Delta(-h,-z_{0})\psi_{V}(v),z_{2}\right)\nonumber\\
& &\cdot \Delta(h,-z_{2})\Delta(-h,-(z_{2}+z_{0}))w\nonumber\\
&=&z_{2}^{-1}\delta\left(\frac{z_{1}-z_{0}}{z_{2}}\right)
E^{-}(-h,z_{1})E^{-}(h,z_{2})z_{0}^{\gamma}\cdot\nonumber\\
& &\cdot Y\left(\psi_{V}\phi_{V} Y(\Delta(-h,z_{0})\phi_{V}^{-1} (u),z_{0})
\Delta(-h,-z_{0})\phi_{V}^{-1} (v),
z_{2}\right)\nonumber\\
& &\cdot \Delta(h,-z_{2})\Delta(-h,-(z_{2}+z_{0}))w
\nonumber\\
&=&z_{2}^{-1}\delta\left(\frac{z_{1}-z_{0}}{z_{2}}\right)
E^{-}(-h,z_{1})E^{-}(-h,z_{2})z_{0}^{\gamma}\cdot\nonumber\\
& &\cdot \psi_{W}\phi_{W}Y(Y(\Delta(-h,z_{0})
\phi_{V}^{-1}(u),z_{0})\Delta(-h,-z_{0})\phi_{V}^{-1}(v),
z_{2})\cdot\nonumber\\
& &\cdot \Delta(-h,z_{2})\Delta(-h, -z_{2}-z_{0})w.
\end{eqnarray}
Suppose $h(0)w=\alpha w$. Then
\begin{eqnarray}
& &z_{2}^{-1}\delta\left(\frac{z_{1}-z_{0}}{z_{2}}\right)
\Delta(-h,-z_{2})\Delta(-h, -z_{2}-z_{0})w\nonumber\\
&=&z_{1}^{-1}\delta\left(\frac{z_{2}+z_{0}}{z_{1}}\right)
(-z_{2}-z_{0})^{-\alpha}\Delta(-h,-z_{2})E^{+}(h,z_{2}+z_{0})w\nonumber\\
&=&z_{1}^{-1}\delta\left(\frac{z_{2}+z_{0}}{z_{1}}\right)
(-z_{2}-z_{0})^{-\alpha}\Delta(-h,-z_{2})E^{+}(h,z_{1})w\nonumber\\
&=&z_{1}^{-1}\delta\left(\frac{z_{2}+z_{0}}{z_{1}}\right)
(-z_{1})^{\alpha}(-z_{2}-z_{0})^{-\alpha}\Delta(-h,-z_{2})\Delta(-h,-z_{1})w\nonumber\\
&=&z_{1}^{-1}\delta\left(\frac{z_{2}+z_{0}}{z_{1}}\right)
\left(\frac{z_{2}+z_{0}}{z_{1}}\right)^{-\alpha}
\Delta(-h,-z_{2})\Delta(-h,-z_{1})w\nonumber\\
&=&z_{2}^{-1}\delta\left(\frac{z_{1}-z_{0}}{z_{2}}\right)
\left(\frac{z_{1}-z_{0}}{z_{2}}\right)^{\alpha}
\Delta(-h,-z_{2})\Delta(-h,-z_{1})w.
\end{eqnarray}
Thus
\begin{eqnarray}
& &z_{2}^{-1}\delta\left(\frac{z_{1}-z_{0}}{z_{2}}\right)
Y(Y(u,z_{0})v,z_{2})w\nonumber\\
&=&z_{2}^{-1}\delta\left(\frac{z_{1}-z_{0}}{z_{2}}\right)
E^{-}(-h,z_{0}+z_{2})E^{-}(-h,z_{2})z_{0}^{\gamma}\cdot\nonumber\\
& &\cdot \psi_{W}\phi_{W}Y(Y(\Delta(-h,z_{0})
\phi_{V}^{-1}(u),z_{0})\Delta(-h,-z_{0})\phi_{V}^{-1}(v), z_{2})
\left(\frac{z_{1}-z_{0}}{z_{2}}\right)^{\alpha}C.
\end{eqnarray}
Thus
\begin{eqnarray}
& &z_{0}^{-1}\delta\left(\frac{z_{1}-z_{2}}{z_{0}}\right)
Y(u,z_{1})Y(v,z_{2})w
-(-1)^{\gamma}z_{0}^{-1}\delta\left(\frac{z_{2}-z_{1}}{-z_{0}}\right)
Y(v,z_{2})Y(u,z_{1})w\nonumber\\
&=&z_{2}^{-1}\delta\left(\frac{z_{1}-z_{0}}{z_{2}}\right)
\left(\frac{z_{1}-z_{0}}{z_{2}}\right)^{\alpha}Y(Y(u,z_{0})v,z_{2})w
\end{eqnarray}
if $h(0)w=\alpha w$.

Case 2: $w\in \tilde{W}$.
Similarly, for $u,v\in \tilde{V}, \tilde{w}\in \tilde{W}$, we have:
\begin{eqnarray}
& & Y(u,z_{1})Y(v,z_{2})\tilde{w}\nonumber\\
&=&E^{-}(-h,z_{1})\phi_{W}Y(\phi_{V}^{-1}u,z_{1})\Delta(-h,-z_{1})\cdot
\nonumber\\
& &\cdot (-1)^{\gamma}E^{-}(-h,z_{2})\psi_{W}Y(
\phi_{V}^{-1}v,z_{2})\Delta(-h,-z_{2})\tilde{w}\nonumber\\
&=&(z_{1}-z_{2})^{\gamma}\phi_{W}\psi_{W}E^{-}(-h,z_{1})E^{-}(-h,z_{2})
Y(\Delta(-h,z_{1}-z_{2})\phi_{V}^{-1}(u),z_{1})\cdot\nonumber\\
& &\cdot Y(\Delta(-h,-z_{1}+z_{2})\phi_{V}^{-1}v,z_{2})\Delta(-h,-z_{1})
\Delta(-h,-z_{2})\tilde{w}.
\end{eqnarray}
Then
\begin{eqnarray}
& &z_{0}^{-1}\delta\left(\frac{z_{1}-z_{2}}{z_{0}}\right)
Y(u,z_{1})Y(v,z_{2})\tilde{w}\nonumber\\
&=&z_{0}^{-1}\delta\left(\frac{z_{1}-z_{2}}{z_{0}}\right)
z_{0}^{\gamma}\phi_{W}\psi_{W}E^{-}(-h,z_{1})E^{-}(-h,z_{2})
Y(\Delta(-h,z_{0})\phi_{V}^{-1}u,z_{1})\cdot\nonumber\\
& &\cdot Y(\Delta(-h,-z_{0})\phi_{V}^{-1}v,z_{2})\Delta(-h,-z_{1})
\Delta(-h,-z_{2})\tilde{w},
\end{eqnarray}
and
\begin{eqnarray}
& &z_{0}^{-1}\delta\left(\frac{z_{2}-z_{1}}{-z_{0}}\right)
Y(v,z_{2})Y(u,z_{1})\tilde{w}\nonumber\\
&=&z_{0}^{-1}\delta\left(\frac{z_{2}-z_{1}}{-z_{0}}\right)
(-z_{0})^{\gamma}\phi_{W}\psi_{W}E^{-}(-h,z_{1})E^{-}(-h,z_{2})
Y(\Delta(-h,-z_{0})\phi_{V}^{-1}v,z_{2})\cdot\nonumber\\
& &\cdot Y(\Delta(-h,-z_{0})\phi_{V}^{-1}u,z_{1})\Delta(-h,-z_{1})
\Delta(-h,-z_{2})\tilde{w}.
\end{eqnarray}
Thus
\begin{eqnarray}
& &\ \ \ z_{0}^{-1}\delta\left(\frac{z_{1}-z_{2}}{z_{0}}\right)
Y(u,z_{1})Y(v,z_{2})\tilde{w}\nonumber\\
& &\ \ \ -(-1)^{\gamma}z_{0}^{-1}\delta\left(\frac{z_{2}-z_{1}}{-z_{0}}\right)
Y(v,z_{2})Y(u,z_{1})\tilde{w}\nonumber\\
& &=z_{2}^{-1}\delta\left(\frac{z_{1}-z_{0}}{z_{2}}\right)
\phi_{W}\psi_{W}E^{-}(-h,z_{1})E^{-}(-h,z_{2})\cdot\nonumber\\
& &\ \ \ \cdot Y\left(Y(\Delta(-h,z_{0})\phi_{V}^{-1}u,z_{0})
\Delta(-h,-z_{0})\phi_{V}^{-1}v,
z_{2}\right)\Delta(-h,-z_{1})
\Delta(-h,-z_{2})\tilde{w}.\;\;\;\;
\end{eqnarray}
On the other hand, using the calculation in case 1, we obtain
\begin{eqnarray}
& &z_{2}^{-1}\delta\left(\frac{z_{1}-z_{0}}{z_{2}}\right)
Y(Y(u,z_{0})v,z_{2})\tilde{w}\nonumber\\
&=&z_{2}^{-1}\delta\left(\frac{z_{1}-z_{0}}{z_{2}}\right)
E^{-}(-h,z_{1})E^{-}(h,z_{2})z_{0}^{\gamma}\cdot\nonumber\\
& &\cdot Y\left(\psi_{V}\phi_{V} Y(\Delta(-h,z_{0})\phi_{V}^{-1}u,z_{0})
\Delta(-h,-z_{0})\phi_{V}^{-1}v,
z_{2}\right)\nonumber\\
& &\cdot \Delta(h,-z_{2})\Delta(-h,-(z_{2}+z_{0}))\tilde{w}
\nonumber\\
&=&z_{2}^{-1}\delta\left(\frac{z_{1}-z_{0}}{z_{2}}\right)
E^{-}(-h,z_{1})E^{-}(h,z_{2})z_{0}^{\gamma}\psi_{W}^{-1}\cdot\nonumber\\
& &\cdot Y\left(\Delta(h,z_{2})\psi_{V}\phi_{V} Y(\Delta(-h,z_{0})
\phi_{V}^{-1}u,z_{0})
\Delta(-h,-z_{0})\phi_{V}^{-1}v,
z_{2}\right)\nonumber\\
& &\cdot \psi_{W}\Delta(h,-z_{2})\Delta(-h,-(z_{2}+z_{0}))\tilde{w}
\nonumber\\
&=&z_{2}^{-1}\delta\left(\frac{z_{1}-z_{0}}{z_{2}}\right)
E^{-}(-h,z_{1})E^{-}(h,z_{2})z_{0}^{\gamma}z_{2}^{-2\gamma}\psi_{W}^{-1}
\cdot\nonumber\\
& &\cdot Y\left(\psi_{V}\phi_{V}\Delta(h,z_{2}) Y(\Delta(-h,z_{0})
\phi_{V}^{-1}u,z_{0})
\Delta(-h,-z_{0})\phi_{V}^{-1}v,
z_{2}\right)\nonumber\\
& &\cdot \psi_{W}\Delta(h,-z_{2})\Delta(-h,-(z_{2}+z_{0}))\tilde{w}
\nonumber\\
&=&z_{2}^{-1}\delta\left(\frac{z_{1}-z_{0}}{z_{2}}\right)
E^{-}(-h,z_{0}+z_{2})E^{-}(-h,z_{2})z_{0}^{\gamma}z_{2}^{-2\gamma}\cdot\nonumber\\
& &\cdot \psi_{W}\phi_{W} Y(\Delta(h,z_{2})Y(\Delta(-h,z_{0})
\phi_{V}^{-1}u,z_{0})\Delta(-h,-z_{0})\phi_{V}^{-1}v, z_{2})
\Delta(-2h,-z_{2})\psi_{W}G\nonumber\\
&=&z_{2}^{-1}\delta\left(\frac{z_{1}-z_{0}}{z_{2}}\right)
E^{-}(-h,z_{0}+z_{2})E^{-}(-h,z_{2})z_{0}^{\gamma}\cdot\nonumber\\
& &\cdot \phi_{W} Y(\Delta(h,z_{2})Y(\Delta(-h,z_{0})
\phi_{V}^{-1}u,z_{0})\Delta(-h,-z_{0})\phi_{V}^{-1}v, z_{2})
\left(\frac{z_{2}+z_{0}}{z_{1}}\right)^{h(0)}\psi_{W}G\nonumber\\
&=&z_{2}^{-1}\delta\left(\frac{z_{1}-z_{0}}{z_{2}}\right)
E^{-}(-h,z_{0}+z_{2})E^{-}(-h,z_{2})z_{0}^{\gamma}\cdot\nonumber\\
& &\cdot \phi_{W}\psi_{W} Y(Y(\Delta(-h,z_{0})
\phi_{V}^{-1}u,z_{0})\Delta(-h,-z_{0})\phi_{V}^{-1}v, z_{2})
\left(\frac{z_{2}+z_{0}}{z_{1}}\right)^{h(0)}G,
\end{eqnarray}
where $G=\Delta(-h,-z_{1})\Delta(-h,-z_{2})\tilde{w}$.
Then the proof is complete. $\;\;\;\;\;\Box$

{\bf Proposition 3.14.} {\it Under the assumption of Theorem 3.13, if
$\tilde{W}$ as a $V$-module is not isomorphic to $W$, then $W\oplus \tilde{W}$
is an irreducible
$\bar{V}$-module, and  if $\tilde{W}$ as a $V$-module is isomorphic to $W$,
then
$W\oplus \tilde{W}$ is a direct sum of two irreducible $\bar{V}$-modules, each
of which
is isomorphic to $W$ as a $V$-module.}

{\bf Proof.} If $\tilde{W}$ as a $V$-module is not isomorphic to $W$, the
irreducibility of
$W\oplus \tilde{W}$ follows from the argument in Remark 3.10. For the rest of
proof, we assume
that $\tilde{W}$ as a $V$-module is isomorphic to $W$. Let $f$ be a
$V$-isomorphism from $W$ onto
$\tilde{W}$. Then $\phi_{W}^{-1}f\psi_{W}f$ is a $V$-automorphisam of $W$.
Since $W$ is an irreducible $V$-module, from  Schur lemma by multiplying a
constant we
 may choose $f$ such that
\begin{eqnarray}
\phi_{W}^{-1}f\psi_{W}f={\rm id}_{W}.
\end{eqnarray}
By definition, for any $u\in \tilde{V}, w\in W, \tilde{w}\in \tilde{W}$, we
have:
\begin{eqnarray}
Y(u,z)w&=&E^{-}(-h,z)\psi_{\tilde{W}}\pi_{W}Y(\phi_{V}^{-1}(u),z)\Delta(-h,-z)w,\\
Y(u,z)\tilde{w}&=&E^{-}(-h,z)\psi_{W}
Y(\phi_{V}^{-1}u,z)\Delta(-h,-z)\tilde{w}.
\end{eqnarray}
Since $f(w)\in \tilde{W}$ by Lemma 3.12 we get
\begin{eqnarray}
& &\ \ \ fY(u,z)f(w)\nonumber\\
& &=E^{-}(-h,z)f\psi_{W}
Y(\phi_{V}^{-1}u,z)\Delta(-h,-z)f(w)\nonumber\\
& &=E^{-}(-h,z)f\psi_{W}
Y(\phi_{V}^{-1}u,z)\Delta(-h,-z)f(w)\nonumber\\
& &=E^{-}(-h,z)f\psi_{W}fY(\phi_{V}^{-1}u,z)(w)\nonumber\\
& &=E^{-}(-h,z)\phi_{W}Y(\phi_{V}^{-1} u,z)\Delta(-h,-z)(w)\nonumber\\
& &=Y(u,z)w.
\end{eqnarray}
Therefore $W_{0}=\{ (w, f(w))|w\in W\}$ and $W_{1}=\{(w,-f(w))|w\in W\}$ are
submodules of
$\bar{W}=W\oplus \tilde{W}$ and $\bar{W}$ is the direct sum of $W_{0}$ and
$W_{1}$. Since
$W_{0}$
 and $W_{1}$ are irreducible $V$-modules, they are irreducible
$\bar{V}$-modules. Then
proof is
complete. $\;\;\;\;\Box$

{\bf Remark 3.15:} Let $V$ be a vertex operator algebra, let $M$ be a
$V$-module and let $U$ be any subspace of $M$. Then it follows from the
associativity
([DL], [FHL]) that the
linear span $V\cdot U$  of all $a_{n}U$ for $a\in V, n\in {\bf Z}$ is a
submodule of
$M$ (cf. Lemma 6.1.1 [Li4]). Furthermore, if $V$ is a simple vertex operator
algebra,
 then  $Y(a,z)u\ne 0$ for $0\ne a\in V, 0\ne u\in M$ ([DL]). This is due to the
fact
that
the annihilating space $Ann_{V}(u)=\{a\in V|Y(a,z)u=0\}$ is an ideal of $V$.

Next we classify all irreducible $\bar{V}$-modules. Recall that $\sigma$ is
an involution of $\bar{V}$. For any $\bar{V}$-module $(M,Y_{M}(\cdot,z))$, we
have a
$\bar{V}$-module $(M,Y_{M}(\sigma\cdot,z))$. For convenience, we denote this
module by
$M^{\sigma}$.

{\bf Lemma 3.16.} {\it If $M$ is a $\bar{V}$-module such that $M$ is an
irreducible
$V$-module,
then $M^{\sigma}$ is not isomorphic to $M$ as a $\bar{V}$-module.}

{\bf Proof.} Suppose $M^{\sigma}$ as a $\bar{V}$-module is isomorphic to $M$,
{\it i.e.,} there
is a linear automorphism $\psi$ of $M$ satisfying
\begin{eqnarray}
\psi(Y_{M}(a,z)u)&=&Y_{M}(a,z)\psi(u)\;\;\;\mbox{for }a\in V,u\in M,\\
\psi(Y_{M}(a,z)u)&=&-Y_{M}(a,z)\psi(u)\;\;\;\mbox{for }a\in \tilde{V},u\in M.
\end{eqnarray}
Since $\psi$ is an automorphism of $M$ as a $V$-module, by Schur lemma
$\psi=\alpha {\rm id}_{M}$ for some nonzero number $\alpha$. Then we have:
\begin{eqnarray}
\alpha Y_{M}(a,z)u=-\alpha Y_{M}(a,z)u\;\;\;\mbox{for }a\in \tilde{V},u\in M.
\end{eqnarray}
Thus $Y_{M}(a,z)u=0$ for any $a\in \tilde{V},u\in M$. It is a
contradiction.$\;\;\;\;\Box$

{\bf Proposition 3.17.}\footnote{It was pointed out by Professor Dong that
it also follows from Theorem 5.4 in [DM].}
 {\it Let $M_{1}$ and $M_{2}$ be $\bar{V}$-modules such that $M_{1}$ and
$M_{2}$ are isomorphic irreducible $V$-modules. Then $M_{1}$ is isomorphic to
either
$M_{2}$ or $M_{2}^{\sigma}$.}

{\bf Proof.} Let $f$ be a $V$-isomorphism from $M_{1}$ onto $M_{2}$. Then
$f^{-1}Y_{2}(\cdot,z)f$ is an intertwining operator of type
$\left(\begin{array}{c}M_{1}\\ \tilde{V} M_{1}\end{array}\right)$.
Since $\tilde{V}$ is a simple current $V$-module, the corresponding fusion rule
is one, so
that there is a nonzero number $\alpha$ such that
\begin{eqnarray}
Y_{1}(a,z)u=\alpha f^{-1}Y_{2}(a,z)f(u)\;\;\;\mbox{ for any }a\in
\tilde{V},u\in M_{1}.
\end{eqnarray}
Then
\begin{eqnarray}
&
&Y_{1}(a,z_{1})Y_{1}(b,z_{2})u=\alpha^{2}f^{-1}Y_{2}(a,z_{1})Y_{2}(b,z_{2})f(u),\\
&
&Y_{1}(b,z_{2})Y_{1}(a,z_{1})u=\alpha^{2}f^{-1}Y_{2}(b,z_{2})Y_{2}(a,z_{1})f(u),\\
& &Y_{1}(Y(a,z_{0})b,z_{2})u=Y_{2}(Y(a,z_{0})b,z_{2})u
\end{eqnarray}
for any $a,b\in \tilde{V},u\in M_{1}$. From the Jacobi identity we obtain
$\alpha^{2}=1$.
Thus $\alpha=1$ or $\alpha=-1$. Therefore $M_{1}$ is isomorphic to either
$M_{2}$ or
$M_{2}^{\sigma}$.$\;\;\;\;\Box$

{\bf Proposition 3.18.} {\it Let $W_{1}=W_{11}\oplus W_{12}$ and
$W_{2}=W_{21}\oplus W_{22}$
 be two $\bar{V}$-modules satisfying the following conditions:
(1) Each $W_{ij}$ is an irreducible $V$-module. (2) $W_{11}$ as a $V$-module is
isomorphic to $W_{21}$. (3) $a_{n}W_{i1}\subseteq W_{i2},\;a_{n}W_{i2}\subseteq
W_{i1}$
for $a\in \tilde{V}, n\in {\bf Z}, i=1,2$. Then $W_{1}$ and $W_{2}$ are
isomorphic $\bar{V}$-modules.}

{\bf Proof.} Let $f_{1}$ be a $V$-isomorphism from $W_{11}$ onto $W_{21}$.
Since
$Y_{2}(\cdot,z)f_{1}$ is an intertwining operator of type
$\left(\begin{array}{c}W_{22}\\ \tilde{V} W_{11}\end{array}\right)$ and
$(W_{12}, Y_{1}(\cdot,z))$ is
a tensor product of $(\tilde{V},W_{11})$, there is a $V$-homomorphism $f_{2}$
from $W_{21}$ to
$W_{22}$ such that
\begin{eqnarray}
Y_{2}(a,z)f_{1}(u)=f_{2}Y_{1}(a,z)u\;\;\;\mbox{for any }a\in \tilde{V}, u\in
W_{11}.
\end{eqnarray}
Let $f=f_{1}\oplus f_{2}$ be the $V$-homomorphism from $W_{1}$ to $W_{2}$. Then
\begin{eqnarray}
f(Y_{1}(a,z)u)&=&Y_{2}(a,z)f(u)\;\;\;\mbox{for any }a\in V_{0}, u\in W_{1},\\
f(Y_{1}(a,z)v)&=&Y_{2}(a,z)f(u)\;\;\;\mbox{for any }a\in V_{1}, v\in W_{11}.
\end{eqnarray}
For any $a,b\in \tilde{V}, u\in W_{11}$, from the associativity (cf. [DL],
[FHL]) we have:
\begin{eqnarray}
f(Y_{1}(a,z_{1})Y_{1}(b,z_{2})u)=Y_{1}(a,z_{1})Y_{1}(b,z_{2})f(u).
\end{eqnarray}
Since $Y_{1}(a,z_{1})Y_{1}(b,z_{2})f(u)=Y_{1}(a,z_{1})f(Y_{1}(b,z_{2})u)$, we
obtain
\begin{eqnarray}
f(Y_{1}(a,z_{1})Y_{1}(b,z_{2})u)=Y_{1}(a,z_{1})f(Y_{1}(b,z_{2})u).
\end{eqnarray}
Since $W_{12}$ is linearly spanned by all the coefficients of $Y_{1}(b,z)u$ for
$b\in V_{1},u\in W_{11}$, we have:
\begin{eqnarray}
f(Y_{1}(a,z)v)=Y_{2}(a,z)f(v)\;\;\;\mbox{for any }a\in \tilde{V},v\in W_{12}.
\end{eqnarray}
Then $f$ is a $\bar{V}$-homomorphism from $W_{1}$ to $W_{2}$. It is clear that
$f$ is an
isomorphism because each $W_{ij}$ is an irreducible $V$-module.$\;\;\;\;\Box$

{\bf Theorem 3.19.} {\it If the vertex operator algebra $V$ is rational,
then $\bar{V}$ is rational and any irreducible $\bar{V}$-module
is isomorphic to one of those obtained in Proposition 3.14.}

{\bf Proof.} Since $V$ is rational, any $\bar{V}$-module $M$ is a completely
reducible
$V$-module. Let $W$ be any irreducible $V$-submodule of $M$. From Remark 3.10,
$\bar{V}\cdot W$
is a $\bar{V}$-submodule of $M$ and $\bar{V}\cdot W=W+V\cdot W$. It follows
from Corollary 2.10
that $(\tilde{W}, \tilde{Y}(\cdot,z))$ is a tensor product of $(W,\tilde{V})$.
Let
$Y^{t}(\cdot,z)$ be the transpose intertwining operator of $Y(\cdot,z)$. Then
$(\tilde{W},\tilde{Y}^{t}(\cdot,z))$ is a tensor product of $(\tilde{V},W)$
 (Lemma 5.1.6 [Li4]).
Since $\tilde{V}$ and $W$ are $V$-modules, $Y(\cdot,z)$ is an intertwining
operator,
 so that
there is  a unique $V$-homomorphism $f_{1}$ from $\tilde{W}$ to $\tilde{V}\cdot
W$
satisfying
\begin{eqnarray}
f_{1}\tilde{Y}^{t}(\tilde{w},z)a=Y(a,z)u\in \tilde{V}\cdot W \;\;\;\mbox{ for }
\tilde{w}\in \tilde{W}, a\in \tilde{V}.
\end{eqnarray}
Consequently, $\tilde{V}\cdot W$ (isomorphic to $\tilde{W}$) is an irreducible
$V$-module.
If $\tilde{V}\cdot W=W$,
then $W$ is an irreducible $\bar{V}$-submodule of $M$. If $\tilde{W}\ne W$,
then
$\bar{V}\cdot W=W\oplus V\cdot W$. By Lemma 3.16, $\bar{V}\cdot W$ is a
completely
reducible $\bar{V}$-module.  Then $M$ is a direct sum of
irreducible $\bar{V}$-modules.$\;\;\;\;\Box$

Let us return to the case for an affine Lie algebra $\tilde{{\bf g}}$.

{\bf Theorem 3.20.} {\it Let $\ell$ be a rational number and let
$\lambda_{i}$ be a cominimal weight such that
$L(\ell,\ell\lambda_{i})$ is self-dual and that
$\gamma=\ell \langle \lambda_{i},\lambda_{i}\rangle$ is an integer.
Then $L(\ell,0)\oplus L(\ell,\ell\lambda_{i})$ is a
 vertex operator algebra if $\gamma$ is even and
$L(\ell,0)\oplus L(\ell,\ell\lambda_{i})$ is a vertex operator
superalgebra if $\gamma$ is odd.}

{\bf Proof.} It is clear that
$L(\ell,0)$ is self-dual for any $\ell$.
Since
$L(\ell,\ell\lambda_{i})$ is self-dual, there is a nonzero intertwining
operator of type
$\left(\!\!\begin{array}{c}L(\ell,0)\\L(\ell,\ell\lambda_{i})
L(\ell,\ell\lambda_{i})\end{array}\!\!\right)$.
Because $L(\ell,\ell\lambda_{i})$ is a simple current, by Corollary 2.10
$\tilde{L}(\ell,\ell\lambda_{i})$
is isomorphic to $L(\ell,0)$. Then it follows from Theorem 3.9
immediately.$\;\;\;\;\Box$

Next we go back to a specific example.
Let $\{e,f,h\}$ be a basis for $sl_{2}$ such that
$[h,e]=e,\;[h,f]=-f,\; [e,f]=h$.
Let $\tilde{sl}_{2}$ be the affine Lie algebra with respect to the
normalized Killing form such that the square length of $h$ is ${1\over 2}$. For
any complex numbers $\ell$ and $j$, denote by $L(\ell,j)$ the irreducible
highest weight module of level $\ell$ with spin $j$ for
$\tilde{sl}_{2}$. It is well
known ([DL], [FZ], [Li2]) that if $\ell$ is a positive integer,
$L(\ell,0)$ is rational and
the set $\{ L(\ell,j)|2j\in {\bf N}, 2j\le \ell\}$ of all
standard $\tilde{sl}_{2}$-modules of level $\ell$ is
the set of equivalence classes of irreducible $L(\ell,0)$-modules.

{\bf Corollary 3.21.}  {\it (a) Let $\ell$ be any complex number such that
$\ell\ne -2$. Then $L(\ell,{\ell\over 2})$ has a natural $L(\ell,0)$-module
structure.}

{\it (b) Let $k$ be a positive integer and let $\ell=4k$.
Then $L(\ell,0)\oplus L(\ell, 2k)$ has a natural vertex operator
algebra structure with an order-two automorphism $\sigma$ such that
$\sigma|_{L(\ell,0)}={\rm id}$, $\sigma|_{L(\ell,2k)}=-{\rm id}$.}

{\bf Proof.} Since $\displaystyle{\langle h,h\rangle ={1\over 2}}$, we get
$\gamma =\ell \langle h,h\rangle ={1\over 2}\ell.$
Then it immediately follows from Theorem 3.9.$\;\;\;\;\Box$

{\bf Corollary 3.22.} {\it Let $k$ be a positive integer and let $\ell=4k$.
Then
(a) the spaces $L(\ell,j)\oplus
L(\ell,2k-j)$ for $j\in {\bf Z},0\le j<k$, $L(\ell,k)$ and $L(\ell,k)^{\sigma}$
have natural irreducible
module structures for the vertex operator algebra $L(\ell,0)\oplus
L(\ell, 2k)$.}

{\it (b) $L(\ell,j)\oplus
L(\ell,2k-j)$ for $j\in {1\over 2}+{\bf Z},0\le j<k$ have
natural irreducible
$\sigma$-twisted module structures for the vertex operator algebra
$L(\ell,0)\oplus L(\ell, 2k)$.}

{\it (c) The vertex
operator algebra $L(\ell,0)\oplus L(\ell, 2k)$ is rational and
any irreducible module for $L(\ell,0)\oplus
L(\ell, 2k)$ is isomorphic to one of the modules in (a).}

{\bf Remark 3.23:} In the physical references (cf. [MSe]), the existence of
twisted modules in Corollary 3.22 (b) was neglected.


\section{Twisted modules for inner automorphisms}
In this section, we present two constructions of twisted modules from any
untwisted
modules for an inner
 automorphism of a vertex operator algebra and we prove that they are
essentially equivalent.

Let $h\in V$ satisfying (\ref{e2.15}).
Furthermore, we assume that $h(0)$ semisimply acts on $V$ with
rational eigenvalues.
It is clear that $e^{2\pi ih(0)}$ is an automorphism of $V$. If the
denominators of all eigenvalues of $h(0)$ are bounded, then
 $e^{2\pi ih(0)}$ is of finite order.

{\bf Proposition 4.1 [Li3].} {\it Let $(M,Y_{M}(\cdot,z))$ be any
$V$-module. Then $(\tilde{M},\tilde{Y}(\cdot,z))=(M, Y(\Delta(h,z)\cdot,z))$ is
a $\sigma$-twisted $V$-module, where $\sigma=e^{2\pi ih(0)}$.}

This is the first construction for twisted modules from any untwisted modules
for
an inner automorphism. This construction is also closely related to shifted
vertex
operators in [Le].

Let $V=\oplus_{n\in {\bf Q}}V_{(n)}$ be a ${\bf Q}$-graded vertex operator
algebra $V$,
i.e., $V$ satisfies all axioms of a
vertex operator algebra except that it may have nonintegral weights.
Let $\sigma=e^{2\pi iL(0)}$. For any
$a\in V_{(n)},n\in {\bf Z}$, we have: $\sigma (a)=a$.
In particular, $\sigma ({\bf 1})={\bf 1}$ and $\sigma (\omega)=\omega$. For
any $a\in V_{(\alpha)},b\in V_{(\beta)}, \alpha,\beta \in {\bf Q}$, we have:
\begin{eqnarray}
e^{2\pi iL(0)}(a_{n}b)&=&e^{2({\rm wt}a+{\rm wt
}b-n-1)\pi i}(a_{n}b)\nonumber\\
&=&e^{2({\rm wt}a+{\rm wt }b)\pi i}(a_{n}b)\nonumber\\
&=&(e^{2\pi iL(0)}a)_{n}(e^{2\pi iL(0)}b).
\end{eqnarray}
Therefore $\sigma=e^{2\pi i L(0)}$ is an
automorphism of $V$ as a vertex operator algebra.

It is clear that $\sigma_{V}$ is of finite order if and only if there is a
positive integer $T$ such that all weights are contained in ${1\over
T}{\bf Z}$. In general, $\sigma$ may be
of infinite order.
Let $M=\oplus_{\alpha\in {\bf C}}M_{(\alpha)}$ be any $V$-module and let
$M'=\oplus_{\alpha\in {\bf C}}M_{(\alpha)}^{*}$ be the restricted dual
[FHL] of $M$. Define
\begin{eqnarray}
\langle Y(a,z)u',v\rangle = \langle u',Y(e^{zL(1)}e^{\pi i L(0)}
z^{-2L(0)}a,z^{-1})v\rangle
\end{eqnarray}
for any $a\in V, u'\in M', v\in M$.

{\bf Proposition 4.2.} {\it Let $V$ be a ${\bf
Q}$-graded vertex operator algebra and let $M$ be a $V$-module. Then
$(M', Y(\cdot,z))$ defined above is a $\sigma^{2}$-twisted $V$-module.}

{\bf Proof.} It is a slight generalization of
contragredient module theory in [FHL]. But for completeness,
we present the details. First we recall the following  formulas from [FHL].
\begin{eqnarray}
z^{L(0)}L(-1)&=&z^{-1}L(-1)z^{L(0)},\\
e^{zL(1)}L(-1)&=&L(-1)e^{zL(1)}+2ze^{zL(1)}L(0)-z^{2}e^{zL(1)}L(1),\\
z^{L(0)}Y(a,z_{0})z^{-L(0)}&=&Y(z^{L(0)}a,z_{0}z),\\
e^{zL(1)}Y(a,z_{0})e^{-zL(1)}&=&Y\left(e^{z(1-zz_{0})}(1-z_{0}z)^{-2L(0)}a,
\frac{z_{0}}{1-zz_{0}}\right).
\end{eqnarray}
Let $a\in V,u'\in M', v\in M$. Then
\begin{eqnarray}
& &{d\over dz}\langle Y(a,z)u',v\rangle\nonumber\\
&=&{d\over dz}\langle
u',Y(e^{zL(1)}e^{\pi iL(0)}z^{-2L(0)}a,z^{-1})v\rangle \nonumber\\
&=&-z^{-2}\langle
u',Y(L(-1)e^{zL(1)}e^{\pi iL(0)}z^{-2L(0)}a,z^{-1})v\rangle
\nonumber\\
& &+\langle u', Y(L(1)e^{zL(1)}e^{\pi iL(0)}z^{-2L(0)}a,z^{-1})v
\rangle\nonumber\\
& &-2z^{-1}\langle u',
Y(e^{zL(1)}e^{\pi iL(0)}L(0)z^{-2L(0)}a,z^{-1})v.
\end{eqnarray}
On the other hand, we have:
\begin{eqnarray}
& &\langle Y(L(-1)a,z)u',v\rangle\nonumber\\
&=&\langle u',Y(e^{zL(1)}e^{\pi iL(0)}z^{-2L(0)}L(-1)a,z^{-1})v
\rangle \nonumber\\
&=&-z^{-2}\langle u',Y(e^{zL(1)}L(-1)e^{\pi iL(0)}z^{-2L(0)}a,
z^{-1})v\rangle \nonumber\\
&=&-z^{-2}\langle
u',Y(L(-1)e^{zL(1)}e^{\pi iL(0)}z^{-2L(0)}a,z^{-1})v
\rangle \nonumber\\
& &-2z^{-1}\langle
u',Y(e^{zL(1)}L(0)e^{\pi iL(0)}z^{-2L(0)}a,z^{-1})v
\rangle \nonumber\\
& &+\langle u',Y(e^{zL(1)}L(1)e^{\pi iL(0)}z^{-2L(0)}a,z^{-1})v
\rangle.
\end{eqnarray}
Thus
\begin{eqnarray}
{d\over dz}\langle Y(a,z)u',v\rangle =\langle Y(L(-1)a,z)u',v\rangle
\;\;\;\mbox{for }a\in V, u'\in M', v\in M.
\end{eqnarray}
For $a,b\in V$, similar to [FHL] we obtain
\begin{eqnarray}
& &z_{0}^{-1}\delta\left(\frac{z_{1}-z_{2}}{z_{0}}\right)
Y\left(v,z_{2}^{-1}\right)
Y\left(u,z_{1}^{-1}\right)
\nonumber\\
& &-z_{0}^{-1}\delta\left(\frac{z_{2}-z_{1}}{-z_{0}}\right)
Y\left(u,z_{1}^{-1}\right)
Y\left(v,z_{2}^{-1}\right)
\nonumber\\
&=&-(z_{1}z_{2})^{-1}\left(-\frac{z_{0}}{z_{1}z_{2}}\right)^{-1}\delta\left(
\frac{-z_{2}^{-1}+z_{1}^{-1}}{-{z_{0}\over z_{1}z_{2}}}\right)
Y\left(v,z_{2}^{-1}\right)Y\left(u,z_{1}^{-1}\right)\nonumber\\
& &+(z_{1}z_{2})^{-1}\left(-\frac{z_{0}}{z_{1}z_{2}}\right)^{-1}\delta\left(
\frac{z_{1}^{-1}-z_{2}^{-1}}{-{z_{0}\over z_{1}z_{2}}}\right)
Y\left(u,z_{1}^{-1}\right)Y\left(v,z_{2}^{-1}\right)\nonumber\\
&=&(z_{1}z_{2})^{-1}z_{2}\delta\left(\frac{z_{2}^{-1}-{z_{0}\over z_{1}z_{2}}}
{z_{1}^{-1}}\right)Y\left(Y\left(u,
-\frac{z_{0}}{z_{1}z_{2}}\right)v,z_{2}^{-1}\right)\nonumber\\
&=&z_{1}^{-1}\delta\left(\frac{z_{2}+z_{0}}{z_{1}}\right)
Y\left(Y\left(u,-\frac{z_{0}}{z_{1}z_{2}}\right)v,z_{2}^{-1}
\right),
\end{eqnarray}
where
\begin{eqnarray}
u=e^{z_{1}L(1)}e^{\pi iL(0)}z_{1}^{-2L(0)}a,\;
v=e^{z_{2}L(1)}e^{\pi iL(0)}z_{2}^{-2L(0)}b.
\end{eqnarray}
Since the right hand side of the desired twisted Jacobi identity should be
\begin{eqnarray}
z_{1}^{-1}\delta\left(\frac{z_{2}+z_{0}}{z_{1}}\right)
\left(\frac{z_{1}-z_{0}}{z_{2}}\right)^{-2r}
Y\left(e^{z_{2}L(1)}e^{\pi iL(0)}z_{2}^{-2L(0)}Y(a,z_{0})b,z_{2}^{-1}
\right),
\end{eqnarray}
where $r$ is the weight of $a$, i.e., $L(0)a=ra$, it is sufficient to prove
\begin{eqnarray}
& &z_{1}^{-1}\delta\left(\frac{z_{2}+z_{0}}{z_{1}}\right)
\left(\frac{z_{1}-z_{0}}{z_{2}}\right)^{-2r}
e^{z_{2}L(1)}e^{\pi iL(0)}z_{2}^{-2L(0)}Y(a,z_{0})\nonumber\\
&=&z_{1}^{-1}\delta\left(\frac{z_{2}+z_{0}}{z_{1}}\right)
Y\left(e^{z_{1}L(1)}e^{\pi iL(0)}z_{1}^{-2L(0)}a,
-\frac{z_{0}}{z_{1}z_{2}}\right)e^{z_{2}L(1)}e^{\pi iL(0)}z_{2}^{-2L(0)},
\end{eqnarray}
or, equivalently to prove
\begin{eqnarray}
& &z_{1}^{-1}\delta\left(\frac{z_{2}+z_{0}}{z_{1}}\right)
\left(\frac{z_{1}-z_{0}}{z_{2}}\right)^{-2r}\cdot\nonumber\\
& &\cdot e^{z_{2}L(1)}e^{\pi iL(0)}z_{2}^{-2L(0)}Y(a,z_{0})z_{2}^{2L(0)}
e^{-\pi iL(0)}e^{-z_{2}L(1)}\nonumber\\
&=&z_{1}^{-1}\delta\left(\frac{z_{2}+z_{0}}{z_{1}}\right)
Y\left(e^{z_{1}L(1)}e^{\pi iL(0)}z_{1}^{-2L(0)}a,
-\frac{z_{0}}{z_{1}z_{2}}\right).
\end{eqnarray}
By the conjugation formulas (4.3)-(4.6), we have:
\begin{eqnarray}
& &e^{z_{2}L(1)}e^{\pi iL(0)}z_{2}^{-2L(0)}Y(a,z_{0})z_{2}^{2L(0)}
e^{-\pi iL(0)}e^{-z_{2}L(1)}\nonumber\\
&=&e^{z_{2}L(1)}e^{\pi iL(0)}Y\left(z_{2}^{-2L(0)}a,-z_{2}^{-2}z_{0}
\right)
e^{-\pi iL(0)}e^{-z_{2}L(1)}\nonumber\\
&=&e^{z_{2}L(1)}Y\left(e^{\pi iL(0)}z_{2}^{-2L(0)}a,-z_{2}^{-2}z_{0}
\right)e^{-z_{2}L(1)}\nonumber\\
&=&Y\left(e^{z_{2}(1+z_{0}z_{2}^{-1})L(1)}(1+z_{0}z_{2}^{-1})^{-2L(0)}
e^{\pi iL(0)}
z_{2}^{-2L(0)}a,-\frac{z_{2}^{-2}z_{0}}{1+z_{0}z_{2}^{-1}}\right)
\nonumber\\
&=&Y\left(e^{(z_{2}+z_{0})L(1)}e^{\pi iL(0)}
(z_{2}+z_{0})^{-2L(0)}a,-\frac{z_{0}}{z_{2}(z_{2}+z_{0})}\right).
\end{eqnarray}
Thus
\begin{eqnarray}
& &z_{1}^{-1}\delta\left(\frac{z_{2}+z_{0}}{z_{1}}\right)
\left(\frac{z_{1}-z_{0}}{z_{2}}\right)^{-2r}\cdot\nonumber\\
& &\cdot e^{z_{2}L(1)}e^{\pi iL(0)}z_{2}^{-2L(0)}Y(a,z_{0})z_{2}^{2L(0)}
e^{-\pi iL(0)}e^{-z_{2}L(1)}\nonumber\\
&=&z_{1}^{-1}\delta\left(\frac{z_{2}+z_{0}}{z_{1}}\right)
\left(\frac{z_{2}+z_{0}}{z_{1}}\right)^{2r}\cdot\nonumber\\
& &\cdot Y\left(e^{(z_{2}+z_{0})L(1)}e^{\pi iL(0)}
(z_{2}+z_{0})^{-2L(0)}a,-\frac{z_{0}}{z_{2}(z_{2}+z_{0})}\right)\nonumber\\
&=&z_{1}^{-1}\delta\left(\frac{z_{2}+z_{0}}{z_{1}}\right)
Y\left(e^{z_{1}L(1)}e^{\pi iL(0)}
z_{1}^{-2L(0)}a,-\frac{z_{0}}{z_{2}z_{1}}\right).
\end{eqnarray}
Then the proof is complete. $\;\;\;\;\Box$

{\bf Remark 4.3.} Let $(V,\omega)$ be a vertex operator algebra and
let $e$ be a Virasoro element of $V$ such that $e(-1)=L(-1)$ and that
$e(0)$ is semisimple on $V$ with rational eigenvalues. Then $(V, e)$ is a
${\bf Q}$-graded vertex (operator) algebra. By Proposition 4.3, we
obtain a $\sigma^{2}$-twisted weak module for $(V,\omega)$, where
$\sigma=e^{2\pi i e(0)}$.

{\bf Remark 4.4.} For Frenkel, Lepowsky and Meurman's Moonshine module
vertex operator algebra $V^{\natural}$ [FLM], many involutions have been
constructed
in [M] as formal exponentials of the weight-zero component $e(0)$ of the vertex
operator
associated to a Virasoro element $e$. Unfortunately, $e(-1)\ne L(-1)$.
It is very interesting to know how we can modify this procedure
to obtain (irreducible) twisted modules for those involutions.

{\bf Remark 4.5.} Let $V$ be a \voa, let $h\in V$ satisfying (\ref{e2.15}) and
set $e=\omega+L(-1)h$. Then  $Y(e,z)=\sum_{n\in{\bf Z}}e(n)z^{-n-2}$ gives  a
representation of the Virasoro algebra on $V$ such that $e(-1)=L(-1)$ (cf.
[DLinM]). In general, $(V,Y(\cdot,z),e)$ is a ${\bf Q}$-graded \voa.
Combining this new Virasoro element with Proposition 4.2, we get the second
construction
for twisted modules from any untwisted modules for an inner
automorphism.

As a matter of fact, the above two constructions are essentially the same.
The following proposition gives the connection between the two constructions of
twisted
modules from untwisted modules for an inner automorphism.

{\bf Proposition 4.6.} {\it Let $V$ be a vertex operator algebra,
let $h\in V$ satisfy the conditions (\ref{e2.15}), let
$(M,Y_{M}(\cdot,z))$ be a $V$-module and let $(M', Y_{M'}(\cdot,z), \omega)$ be
the
contragredient module of $M$
with respect to the Virasoro element $\omega$. Suppose that the restricted dual
spaces of
$(M, Y_{M}(\cdot,z), \omega)$ and $(M,Y_{M}(\cdot,z), \omega+L(-1)h)$ are the
same.
Then $(M, Y(\Delta(-2h,z)e^{-\pi ih(0)}\cdot,z))$ as a $e^{-4\pi
ih(0)}$-twisted $V$-module
is isomorphic to the contragredient module of
$(M',Y_{M'}(\cdot,z), \omega)$ with respect to the Virasoro element $\omega +
h(-2){\bf 1}$.}

{\bf Proof.} First from [FHL] we have the following formula:
\begin{eqnarray}
z^{L(0)}e^{z_{0}L(1)}=e^{z_{0}z^{-1}L(1)}z^{L(0)}.
\end{eqnarray}
Set $e=\omega +h(-2){\bf 1}$. Then
\begin{eqnarray}
e(m)=L(m)-(m+1)h(m)\;\;\;\mbox{ for }m\in {\bf Z}.
\end{eqnarray}
In particular, we have:
\begin{eqnarray}
e(-1)=L(-1),\; e(0)=L(0)-h(0),\;e(1)=L(1)-h(1).
\end{eqnarray}
For any $a\in V, u\in M, v'\in M'$, we have:
\begin{eqnarray}
& &\langle Y(a,z)u,v'\rangle\nonumber\\
&=&\langle u, Y(e^{zL(1)}(-z^{-2})^{L(0)}a,z^{-1})v'\rangle\nonumber\\
&=&\langle Y(e^{z^{-1}e(1)}(-z^{2})^{e(0)}
e^{zL(1)}(-z^{-2})^{L(0)}a,z)u,v'\rangle.
\end{eqnarray}
Then
\begin{eqnarray}
& &e^{z^{-1}e(1)}(-z^{2})^{e(0)}
e^{zL(1)}(-z^{-2})^{L(0)}\nonumber\\
&=&e^{z^{-1}(L(1)-2h(1))}(-z^{2})^{L(0)-h(0)}
e^{zL(1)}(-z^{-2})^{L(0)}\nonumber\\
&=&z^{-2h(0)}e^{z^{-1}(L(1)-2h(1))}(-z^{2})^{L(0)}
e^{zL(1)}(-z^{-2})^{L(0)}\nonumber\\
&=&z^{-2h(0)}e^{z^{-1}(L(1)-2h(1))}e^{-z^{-1}L(1)}e^{-\pi ih(0)}.
\end{eqnarray}
For $a\in V, u\in M, v'\in M'$,  using Lemma 3.3 we obtain
\begin{eqnarray}
& &\langle Y(a,z)u,v'\rangle\nonumber\\
&=&\langle Y\left(z^{-2h(0)}e^{z^{-1}(L(1)-2h(1))}e^{-z^{-1}L(1)}
e^{-\pi ih(0)}a,z\right)u,v'\rangle \nonumber\\
&=&\langle Y\left(z^{-2h(0)}\exp\left(\sum_{k=1}^{\infty}
\frac{2h(k)}{k}(-z^{-1})^{k}\right)e^{-\pi ih(0)}a,z\right)u,v'\rangle
\nonumber\\
&=&\langle Y\left(\Delta(-2h,z)e^{-\pi ih(0)}a,z\right)u,v'\rangle.
\end{eqnarray}
Then the proof is complete.$\;\;\;\;\Box$

Combining  Proposition 2.16 with Proposition 4.6 we get

{\bf Corollary 4.7.} {\it Let $L$ be a positive-definite even lattice. Then
every
irreducible $V_{L}$-module can be
considered as a contragredient module of the adjoint module with
respect to some Virasoro element.}

\end{document}